\RequirePackage{ifpdf}
\ifpdf 
\documentclass[pdftex]{sigma}
\else
\documentclass{sigma}
\fi

\DeclareMathOperator{\rank}{rank}
\DeclareMathOperator{\Ima}{Im}
\DeclareMathOperator{\Ker}{Ker}
\DeclareMathOperator{\diag}{diag}
\DeclareMathOperator{\lk}{lk}
\DeclareMathOperator{\St}{star}
\DeclareMathOperator{\const}{const}

{\theoremstyle{definition}
\newtheorem{convention}{Convention}}
\begin{document}

\allowdisplaybreaks

\renewcommand{\PaperNumber}{032}

\FirstPageHeading

\ShortArticleName{A Euclidean Geometric Invariant of Framed (Un)Knots in Manifolds}

\ArticleName{A Euclidean Geometric Invariant\\ of Framed (Un)Knots in Manifolds}

\Author{J\'er\^ome DUBOIS~$^\dag$, Igor G. KOREPANOV~$^\ddag$ and Evgeniy V. MARTYUSHEV~$^\ddag$}

\AuthorNameForHeading{J.~Dubois, I.G.~Korepanov and E.V.~Martyushev}

\Address{$^\dag$~Institut de Math\'ematiques de Jussieu,
Universit\'e Paris Diderot--Paris 7,\\
\hphantom{$^\dag$}~UFR de Math\'ematiques,
Case 7012, B\^atiment Chevaleret,
2, place Jussieu,\\
\hphantom{$^\dag$}~75205 Paris Cedex 13, France}
\EmailD{\href{mailto:dubois@math.jussieu.fr}{dubois@math.jussieu.fr}}

\Address{$^\ddag$~South Ural State University, 76 Lenin Avenue, Chelyabinsk 454080, Russia}
\EmailD{\href{mailto:kig@susu.ac.ru}{kig@susu.ac.ru}, \href{mailto:mev@susu.ac.ru}{mev@susu.ac.ru}}

\ArticleDates{Received October 09, 2009, in f\/inal form April 07, 2010;  Published online April 15, 2010}

\Abstract{We present an invariant of a three-dimensional manifold with a framed knot in it based on the Reidemeister torsion of an acyclic complex of Euclidean geometric origin. To show its nontriviality, we calculate the invariant for some framed (un)knots in lens spaces. Our invariant is related to a f\/inite-dimensional fermionic topological quantum f\/ield theory.}


\Classification{57M27; 57Q10; 57R56}

\section{Introduction}

In this paper, a construction of invariant of three-dimensional manifolds with triangulated boundary is presented, on the example of the complement of a tubular neighborhood of a knot in a closed manifold; the boundary triangulation corresponds in a canonical way to a \emph{framing} of the knot. Algebraically, our invariant is based, f\/irst, on some striking dif\/ferential formulas (see \eqref{eq QRQR} and~\eqref{eq QRQR minors} below) corresponding naturally to Pachner moves~-- elementary rebuildings of a manifold triangulation. These formulas involve some geometric values put in correspondence to triangulation simplexes; specif\/ically, we introduce Euclidean geometry in every tetrahedron. Second, it turns out that the relevant context where these formulas work is the theory of \emph{Reidemeister torsions}.

Recall that Reidemeister torsion made its f\/irst appearance in 1935, in the work of Reidemeister~\cite{Reidemeister} on  the combinatorial classif\/ication of the three-dimensional lens spaces by means of the based simplicial chain complex of the universal cover. Our theory, which stems from the discovery of a ``Euclidean geometric'' invariant of three-dimensional manifolds in paper~\cite{JNMP1}, is, however, radically new, since it unites the algebraic construction of Reidemeister torsion with simplex geometrization. Historically, geometrization came f\/irst, and the Reidemeister torsion came into play only in papers~\cite{TMP24,TMP15} (these papers deal mainly with the case of \emph{four}-dimensional manifolds, so it was this more complicated case that pressed us to clarify the algebraic nature of our constructions).

Our invariant was initially proposed in~\cite{JNMP1} for \emph{closed} manifolds. The next natural step is the investigation of these invariants for manifolds \emph{with boundary}. In doing so, we are guided by the idea of constructing eventually a topological quantum f\/ield theory (TQFT) according to some version of Atiyah's axioms~\cite{A} where, as is known, the boundary plays a fundamental role. Some fragments of this theory have been already developed in our works; in particular, it was shown in papers~\cite{korepanov-2009-1} and~\cite{korepanov-2009-2} that, indeed, a TQFT is obtained this way. This TQFT is \emph{fermionic}: the necessary modif\/ication of Atiyah's axioms is that the usual composition of tensor quantities, corresponding to the gluing of manifolds, is replaced with Berezin (fermionic) integral in anti-commuting variables. What is lacking in papers~\cite{korepanov-2009-1} and~\cite{korepanov-2009-2} is, f\/irst, a systematic exposition of the foundations of the theory and, second, interesting nontrivial examples.

This determines the aims of the present paper: we concentrate on a detailed exposition of the foundations, and we restrict ourselves to the simplest nontrivial case: three-dimensional manifold with toric boundary and a specif\/ic triangulation given on this boundary which f\/ixes the meridian and parallel of the torus. Thus, the whole picture can be imagined as a closed manifold with the withdrawn tubular neighborhood of a knot and, moreover, the framing of the knot is f\/ixed. To make the paper self-consistent, we also concentrate on just one invariant~-- it is called ``zeroth-level invariant'' in~\cite{korepanov-2009-1} and~\cite{korepanov-2009-2}. Even this one invariant turns out to be interesting enough; as for the other invariants, forming together a multi-component object needed for a~TQFT, the reader can consult the two mentioned papers to learn how to construct them. Note that the theory in~\cite{korepanov-2009-2} works for a three-manifold with \emph{any} number of boundary components, each being a two-sphere with any number of handles.

The key object in our theory is the matrix $(\partial\omega_a/\partial l_b)$ of partial derivatives of the so-called def\/icit angles~$\omega_a$ with respect to the edge lengths~$l_b$, where subscripts $a$ and~$b$ parametrize the edges (see Section~\ref{sec_deficit_angles} for detailed def\/initions). The invariant considered in~\cite{JNMP1} makes use of the largest nonvanishing minor of matrix $(\partial\omega_a/\partial l_b)$; some special construction was used to eliminate the non-uniqueness in the choice of this minor, and it has been shown later in~\cite[Section~2]{TMP24} and~\cite[Section~2]{TMP15} that this construction consisted, essentially, in taking the torsion of an acyclic complex built from dif\/ferentials of geometric quantities.

There exists also a version of this invariant using the universal cover of the considered manifold~$M$ and nontrivial representations of the fundamental group $\pi_1(M)$ into the group of motions of three-dimensional Euclidean space~\cite{KM}. In this way, an invariant which seems to be related to the usual Reidemeister torsion has been obtained. A good illustration is the following formula for the invariant of lens spaces proved in~\cite{M-thesis}:
\begin{equation}\label{IK}
\mathrm{Inv}_k\left(L(p, q)\right) = - \frac{1}{p^2} \left(4 \sin \frac{\pi k}{p} \sin \frac{\pi k q}{p}\right)^4.
\end{equation}
Here $L(p, q)$ is a three-dimensional lens space; the subscript~$k$ takes integer values from 1 to the integral part of $p/2$; the invariant consists of real numbers corresponding to each of these~$k$. One can check that the expression in parentheses in formula~(\ref{IK}) is, up to a constant factor, nothing but the usual Reidemeister torsion of $L(p,q)$, see, e.g.,~\cite[Theorem~10.6]{Turaev:2001}.

The invariants appearing from nontrivial representations of $\pi_1(M)$ form an important area of research. This applies to the usual Reidemeister torsion for manifolds and knots~\cite{JDFourier} as well as our ``geometric'' torsion. One can f\/ind some conjectures, concerning the relation of ``geometric'' and ``usual'' invariants constructed using Reidemeister torsions and based on computer calculations, in paper~\cite{M}. Note however that the important feature of the present paper is that we are \emph{not} using any nontrivial representation of the manifold fundamental group or knot group. Formula~(\ref{IK}) has been cited here only to illustrate the fact that, in some situations, the invariant obtained from ``geometric'' torsion can be expressed through the usual Reidemeister torsion.

Returning to the present paper, we introduce here, as we have already said, an invariant of a pair consisting of a manifold and a framed knot in it, and then show its nontriviality on some simplest examples of ``unknots'', i.e.\ simplest closed contours, in the sphere~$S^3$ and lens spaces.

Our geometrization of triangulation simplexes is basically the same as in the Regge theory of discrete general relativity~\cite{Regge}. In this connection, we would like to remark that our theory (\emph{unlike}, for instance, the Ponzano--Regge model~\cite{PR})
\begin{itemize}
\itemsep=0pt
	\item is perfectly f\/inite-dimensional (not involving such things as functional integrals or spins taking inf\/inite number of values) and mathematically strict,
	\item admits generalizations to other geometries, involving even \emph{non-metric} ones, see, e.g.,~\cite{SL2,KKM},
	\item can be generalized to manifolds of more than three dimensions, as papers~\cite{TMP33,TMP24,TMP15}, where a~similar Euclidean geometric invariant is constructed in the four-dimensional case, strongly suggest.
\end{itemize}

\subsection*{Organization}

The main part of the paper starts with the technical Section~\ref{sec_pachner}: we show that \emph{relative Pachner moves}~-- those not involving the boundary~-- are enough to come from any triangulation within the manifold to any other one. Thus, any value invariant under these relative moves is an invariant of the manifold with f\/ixed boundary triangulation. In Section~\ref{sec_deficit_angles}, we def\/ine geometric values needed for the construction of an acyclic complex, and in Section~\ref{sec_acyclic_general} we show how to construct this complex and prove the invariance of its Reidemeister torsion, multiplied by some geometric values, with respect to relative Pachner moves. In Section~\ref{sec_framing}, we show how to change the knot framing within our construction, and how this af\/fects the acyclic complex. In the next two sections we consider our examples: framed unknot in three-sphere (Section~\ref{sec_unknot}) and framed ``unknots'' in lens spaces (Section~\ref{sec_lenses}). In Section~\ref{sec_discussion}, we discuss the results of our paper.

\section{Triangulation for a manifold with a framed knot in it\\ and relative Pachner moves}
\label{sec_pachner}

\begin{figure}[ht]
\centerline{\includegraphics[scale=0.3]{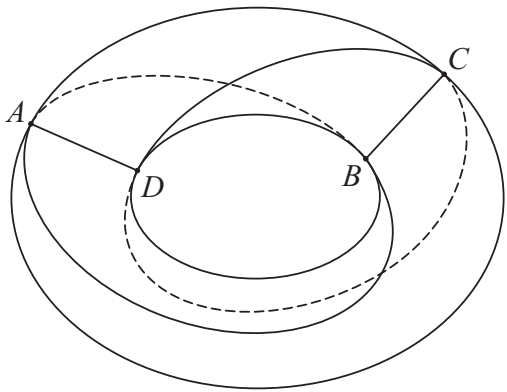}}
\caption{A chain of two identically oriented tetrahedra~$ABCD$.}
\label{fig1a}
\end{figure}

\begin{figure}[ht]
\centerline{\includegraphics[scale=0.3]{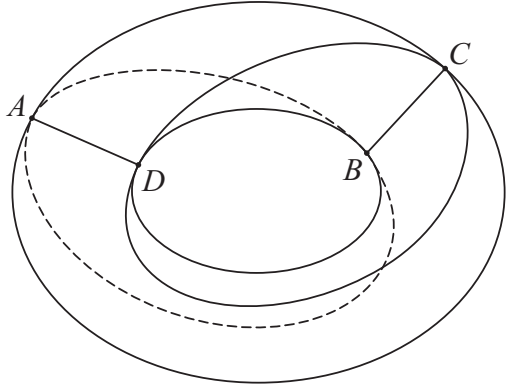}}
\caption{A chain of two oppositely oriented tetrahedra~$ABCD$.}
\label{fig1b}
\end{figure}

We consider a closed oriented three-manifold~$M$ and a triangulation~$T$ of it containing a distinguished chain of two tetrahedra~$ABCD$ of one of the forms depicted in Figs.~\ref{fig1a} and~\ref{fig1b}. These two tetrahedra can either have the same orientation, as in Fig.~\ref{fig1a}, or the opposite orientations, as in Fig.~\ref{fig1b}.

Our construction of the invariant requires adopting the following convention (see Subsection~\ref{subs_daav} for details).

\begin{convention}
\label{conv1}
Any triangulation considered in this paper, including those which appear below at any step of a sequence of Pachner moves, is required to possess the following property: \emph{all} vertices of any tetrahedron are dif\/ferent.
\end{convention}

\begin{remark}
\label{remark-about-convention}
In particular, this convention is satisf\/ied by a \emph{combinatorial triangulation}, i.e., such a triangulation where any simplex is uniquely determined by the set of its vertices, all of those being dif\/ferent.
\end{remark}

To select a special chain of two tetrahedra as depicted in Figs.~\ref{fig1a} and~\ref{fig1b} essentially means the same as to select a \emph{framed knot} in~$M$. To be exact, there is a knot with \emph{two} framings given either by two closed lines (which we imagine as close to each other) $ACA$ and $DBD$, or by the two lines $ABA$ and $DCD$. In the case of the same orientation of the two tetrahedra, these possibilities lead to framings which dif\/fer in one full revolution (of the ribbon between two lines), so, to ensure the invariant character of choosing the framing, we have to choose the ``intermediate'' framing dif\/fering from them both in one-half of a revolution as the framing corresponding to our picture. In the case of the opposite orientations of the two tetrahedra, both ways simply give the same framing.

\begin{remark}
\label{remark-about-half-integers}
Thus, a \emph{half-integer} framing corresponds to Fig.~\ref{fig1a}, represented by a ribbon going like a M\"obius band, and an integer framing~-- to Fig.~\ref{fig1b}. These are the two kinds of framings we will be dealing with in this paper.
\end{remark}

Our aim is to construct an invariant of a pair $(M, K)$, where $K$ is a framed knot in~$M$, starting from a triangulation of~$M$ containing two distinguished tetrahedra as in Fig.~\ref{fig1a} or~\ref{fig1b}. To achieve this, we will construct in Section~\ref{sec_acyclic_general} a value not changing under Pachner moves on triangulation of~$M$ \emph{not touching the distinguished tetrahedra} of Fig.~\ref{fig1a} or~\ref{fig1b}. By ``not touching'' we understand those moves that do not replace either of the two tetrahedra in Fig.~\ref{fig1a} or~\ref{fig1b} with any other tetrahedra, and we call such moves \emph{relative Pachner moves}.

\begin{figure}[ht]
\centerline{\includegraphics[scale=0.39]{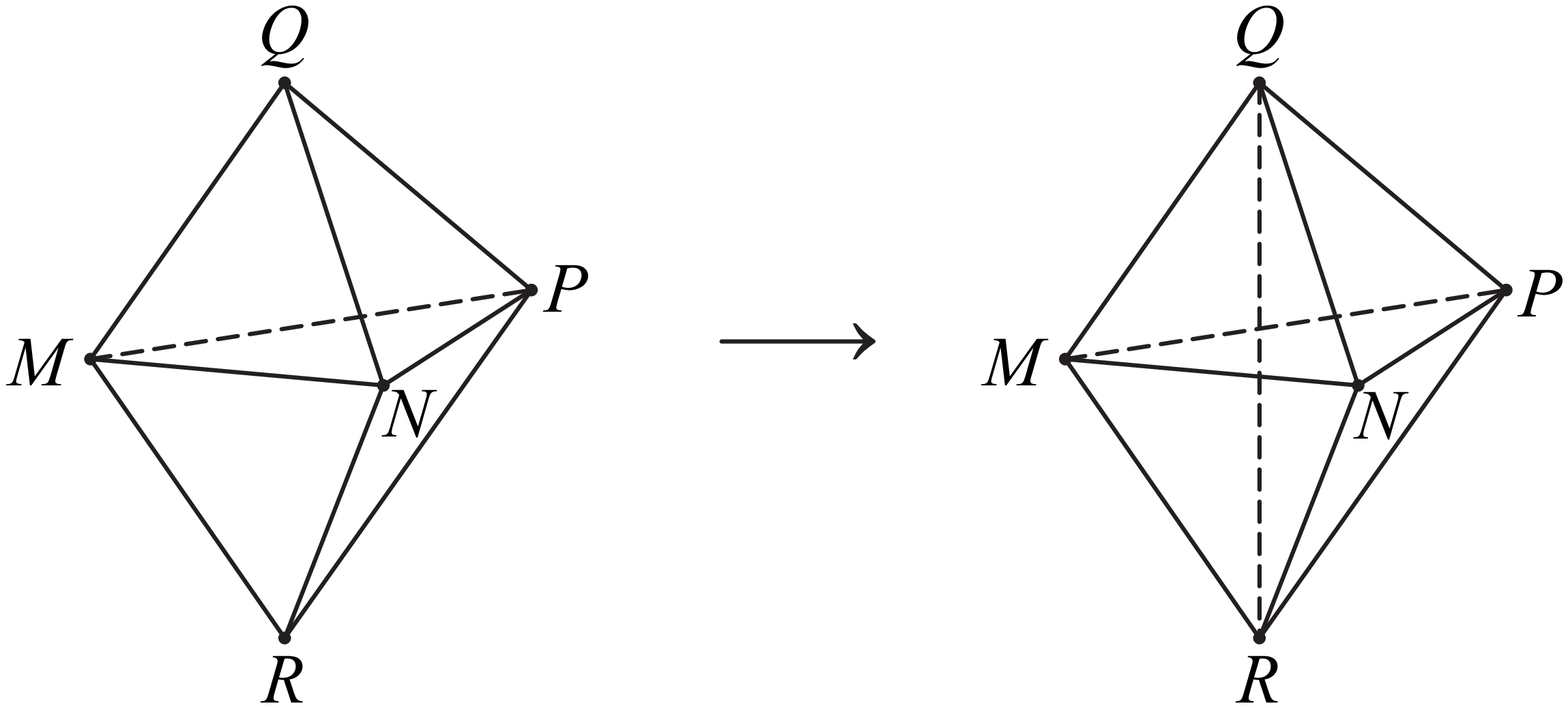}}
\caption{A $2 \to 3$ Pachner move in three dimensions.}
\label{fig2a}
\end{figure}

\begin{figure}[ht]
\centerline{\includegraphics[scale=0.39]{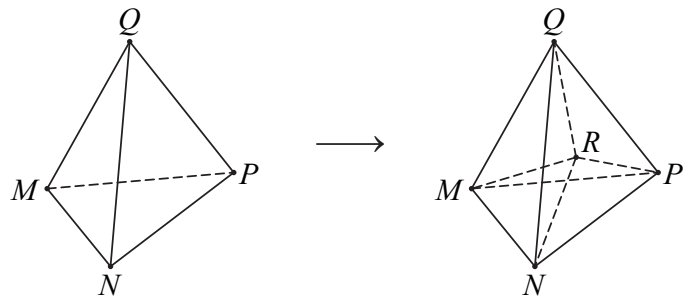}}
\caption{A $1 \to 4$ Pachner move in three dimensions.}
\label{fig2b}
\end{figure}

Recall that Pachner moves are elementary rebuildings of a \emph{closed} triangulated manifold. There are four such moves on three-dimensional manifolds. Two of them are illustrated in Figs.~\ref{fig2a} and~\ref{fig2b}, and the other two are inverse to these. The move in Fig.~\ref{fig2a} replaces the two adjacent tetrahedra $MNPQ$ and $RMNP$ with three new tetrahedra: $MNRQ$, $NPRQ$, and~$PMRQ$. The move in Fig.~\ref{fig2b} replaces one tetrahedron $MNPQ$ with four of them: $MNPR$, $MNRQ$, $MPQR$, and~$NPRQ$.

The main objective of the present section is to prove the following
\begin{theorem}
\label{th1}
Let M be a closed oriented three-manifold, $T_1$ and $T_2$ its triangulations with the same chain of two distinguished tetrahedra~$ABCD$, as depicted in Figs.~{\rm \ref{fig1a}} and~{\rm \ref{fig1b}}. Then, $T_1$~and~$T_2$ are related by a sequence of relative Pachner moves.
\end{theorem}

\begin{proof}
We are going to apply techniques from Lickorish's paper~\cite{L}. Therefore, let us f\/irst explain a method to subdivide the triangulations $T_1$ and $T_2$ in such way that they become combinatorial. Together with Pachner moves (see Figs.~\ref{fig2a} and~\ref{fig2b}), we will use \emph{stellar moves}, see~\cite[Section 3]{L}. In three dimensions, there is no dif\/f\/iculty to express the latter in terms of the former (and vice versa).

We can assume that $T_1$ and~$T_2$ already do not contain any more edges or two-dimensional faces whose vertices all lie in the set $\{A,B,C,D\}$ except those depicted in Figs.~\ref{fig1a} and~\ref{fig1b}~-- otherwise, we can always make obvious stellar subdivisions to ensure this.

Starting from the triangulation $T_1$, we f\/irst do Pachner moves $1\to 4$ in all tetrahedra adjacent to those two in Figs.~\ref{fig1a} and~\ref{fig1b}. Thus, there have appeared eight new vertices, we call them $N_1,\ldots,N_8$.

Next, we look at the edges in Figs.~\ref{fig1a} and~\ref{fig1b}. We are going to make some moves so that the link of each of them contain exactly one vertex between any two~$N_i$. If there were more such vertices, we can eliminate them from the link by doing suitable $2\to 3$ Pachner moves. Namely, the $2\to 3$ Pachner move provides a new edge joining $N_i$ directly with a ``farther'' vertex in the link, thus eliminating from the link the vertex next to $N_i$, see Fig.~\ref{fig3}.
A special case is two edges $AD$ and~$BC$: they require such procedure to be applied twice, ``on two sides''.
\begin{figure}[ht]
\centerline{\includegraphics[scale=0.29]{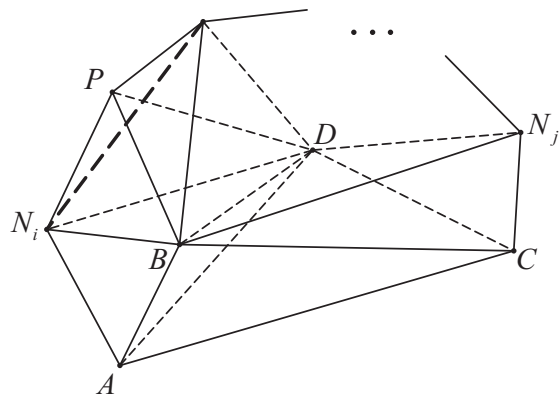}}
\caption{The edge drawn in boldface dashed line appears as a result of move $2\to 3$ and eliminates vertex~$P$ from the link of~$BD$.}
\label{fig3}\vspace{-2mm}
\end{figure}

This done, we make stellar subdivisions in the two-dimensional faces which are joins of the edges in Fig.~\ref{fig1a} or~\ref{fig1b} and the vertices lying between the $N_i$'s~-- two such vertices for each of $AD$ and~$BC$, and one for each of the remaining edges.

After that, we remove from the resulting simplicial complex those tetrahedra that have \emph{at least two} vertices in the set $\{A,B,C,D\}$. Specif\/ically, take these tetrahedra together with all their faces and denote by $L$ the obtained subcomplex. Then we remove~$L$ from our simplicial complex and take the closure of what remains; let~$V_1$ denote this closure.

Note that the triangulated boundary of~$V_1$ (which is also the boundary of~$L$) can be described as follows: f\/irst, \emph{double} the edges $AD$ and $BC$ in Figs.~\ref{fig1a} and~\ref{fig1b} in such way as to make a torus out of the boundary of the tetrahedron chain, and then make a barycentric subdivision of this triangulated torus.

Finally, we subdivide our simplicial complex~$V_1$, doing, e.g., suitable stellar moves in its simplices but leaving the boundary untouched, so that it becomes a combinatorial triangulation, as it is required in order to apply the techniques from~\cite{L}. Let~$W_1$ denote the resulting simplicial complex.

Now, we apply the above procedure to the triangulation $T_2$ and similarly obtain the simplicial complex $W_2$.
Obviously, $W_1$ and $W_2$ are PL-homeomorphic. Then, according to \cite[Theorem~4.5]{L}, these simplicial complexes are stellar equivalent. Moreover, there exists a chain of stellar and inverse stellar moves transforming $W_1$ into $W_2$ and such that \emph{the initial edges in the boundary~$\partial W_1$ never dissapear during the whole process}: they may be at most divided by several starrings done at them, but f\/inally the obtained fragments are glued together again; accordingly, the initial vertices in~$\partial W_1$~-- ends of initial edges~-- remain intact.

Indeed, the cited theorem in~\cite{L} is valid for simplicial complexes and not necessarily mani\-folds. Thus, we can do the following trick: glue to any edge in~$\partial W_1$ an additional two-cell~-- a triangle~-- by one of its sides, and do the same with edges in~$\partial W_2$. The obtained simplicial complexes $W'_1$ and~$W'_2$ are still~PL~-- and consequently stellar~-- homeomorphic, and obviously the additional cells are conserved (at most divided in parts but then glued again) along the whole chain of starrings and inverse starrings, marking thus the edges.

Such chain of starrings and inverse starrings can be extended to the union $W_1\cup L$, so that the result is $W_2\cup L$. Indeed, every stellar move involving~$\partial L$ is done either on one of initial edges or inside one of initial triangles. It is not hard to see that the extension from~$\partial L$ to~$L$ of a starring or inverse starring goes smoothly in both cases.

The subcomplex in Figs.~\ref{fig1a} and~\ref{fig1b} will not be touched by any move in the sequence. So, what remains is to replace the stellar moves with suitable (sequences of) Pachner moves.
\end{proof}

\section{Geometric values needed for the acyclic complex}
\label{sec_deficit_angles}

We are now going to construct an acyclic complex which produces the invariant of a three-manifold with a framed knot in it given by a chain of two tetrahedra as in Figs.~\ref{fig1a} and~\ref{fig1b}. The complex will be like those in \cite{TMP24, TMP15}, but in fact a bit simpler.

\begin{convention}
\label{conv2}
Recall that we are considering an \emph{orientable} manifold~$M$. From now on, we f\/ix a \emph{consistent orientation} for all tetrahedra in the triangulation. The orientation of a tetrahedron is understood here as an ordering of its vertices up to an even permutation; for instance, two tetrahedra $MNPQ$ and~$RMNP$, having a common face~$MNP$, are consistently oriented.
\end{convention}

\subsection{Oriented volumes and def\/icit angles}
\label{subs_daav}
We need the so-called \emph{deficit angles} corresponding to the edges of triangulation. The rest of this section is devoted to explaining these def\/icit angles and related notions, while the acyclic complex itself will be presented in Section~\ref{sec_acyclic_general}.

Recall that we assume that \emph{all} the vertices of any tetrahedron in the triangulation are dif\/ferent (see Convention~\ref{conv1}). Put all the vertices of the triangulation in the Euclidean space $\mathbb R^3$ (i.e., we assign to each of them three real coordinates) in an arbitrary way with only one condition: the vertices of each tetrahedron in the triangulation must not lie in the same plane. This condition ensures that the geometric quantities we will need~-- edge lengths and tetrahedron volumes~-- never vanish.

When we put an oriented tetrahedron $MNRQ$ into~$\mathbb R^3$, we can ascribe to it an \emph{oriented volume} denoted $V_{MNRQ}$ according to the formula
\begin{equation}
6V_{MNRQ}=\overrightarrow{MN} \cdot \overrightarrow{MQ} \cdot \overrightarrow{MR}
\label{eq_oriented_volume}
\end{equation}
(scalar triple product in the right-hand side). If the sign of the volume def\/ined by~(\ref{eq_oriented_volume}) of a given tetrahedron is positive, we say that it is put in~$\mathbb R^3$ \emph{with its orientation preserved}; if it is negative we say that it is put in~$\mathbb R^3$ \emph{with its orientation changed}.

Now we consider the dihedral angles at the edges of triangulation. We will ascribe a sign to each of these angles \emph{coinciding with
the oriented volume sign} of the tetrahedron to which the angle belongs. Consider a certain edge $QR$ in the triangulation, and let its link contain vertices $P_1, \ldots, P_n$, so that the tetrahedra $P_1P_2RQ, \ldots, P_nP_1RQ$ are situated around~$QR$ and form its \emph{star}. With our def\/inition for the signs of dihedral angles, one can observe that the algebraic sum of all angles at the edge~$QR$ is a multiple of $2\pi$, if these angles are calculated according to the usual formulas of Euclidean geometry, starting from \emph{given coordinates of vertices} $P_1, \ldots, P_n$, $Q$ and~$R$.

Namely, we would like to use the following method for computing these dihedral angles. Given the coordinates of vertices, we calculate all the \emph{edge lengths} in tetrahedra $P_1P_2RQ, \ldots$, $P_nP_1RQ$ and the signs of all tetrahedron volumes, and then we calculate dihedral angles from the edge lengths.

Suppose now that we have slightly, but otherwise arbitrarily, changed the edge lengths. Each separate tetrahedron $P_1P_2RQ, \ldots, P_nP_1RQ$ remains still a Euclidean tetrahedron, but the algebraic sum of their dihedral angles at the edge~$QR$ ceases, generally speaking, to be a~multiple of $2\pi$. This means that these tetrahedra can no longer be put in~$\mathbb R^3$ together. In such situation, we call this algebraic sum, taken with the opposite sign, \emph{deficit (or defect) angle} at edge~$QR$:
\begin{equation}
\omega_{QR} = - \sum_{i=1}^n \varphi_i \bmod 2\pi,
\label{eq_deficit_angle}
\end{equation}
where $\varphi_i$ are the dihedral angles at $QR$ in the $n$ tetrahedra under consideration. Note that the minus sign in~(\ref{eq_deficit_angle}) is just due to a convention in ``Regge calculus''~\cite{Regge} where such def\/icit angles often appear.

\subsection{A relation for inf\/initesimal deformations of def\/icit angles}

To build our chain complex~(\ref{eq_alg_complex}) in Section~\ref{sec_acyclic_general}, we need only \emph{infinitesimal} def\/icit angles arising from inf\/initesimal deformations of edge lengths in the neighborhood of a f\/lat case, where all def\/icit angles vanish.

\begin{lemma}
\label{lemma_infinitesimal}
Let $Q$ be any vertex in a triangulation of closed oriented three-manifold, and $QR_1,\ldots,\allowbreak QR_m$ all the edges that end in~$Q$. If infinitesimal deficit angles $d\omega_{QR_i}$ are obtained from infinitesimal deformations of length of edges in the triangulation with respect to the flat case, then{\samepage
\begin{equation}
\sum_{i=1}^m \vec e_{QR_i}\, d\omega_{QR_i} = \vec 0,
\label{eq_zero_for_dx*}
\end{equation}
where $\vec e_{QR_i}=1/l_{QR_i} \cdot \overrightarrow{QR_i}$ is a unit vector along $\overrightarrow{QR_i}$.}
\end{lemma}

\begin{proof}
We f\/irst consider the case of small but f\/inite deformations of edge lengths and def\/icit angles.

Let the edge lengths in tetrahedra $P_1P_2R_iQ, \ldots, P_nP_1R_iQ$, forming the star of edge~$QR_i$, be slightly deformed with respect to the f\/lat case $\omega_{QR_i}=0$. We can introduce a Euclidean coordinate system in tetrahedron~$P_1P_2R_iQ$. Then, this coordinate system can be extended to tetrahedron~$P_2P_3R_iQ$ through their common face~$P_2R_iQ$. Continuing this way, we can go around the edge~$QR_i$ and return in the initial tetrahedron~$P_1P_2R_iQ$, obtaining thus a new coordinate system in it. The transformation from the old system to the new one is an orthogonal rotation around the edge~$QR_i$ through the angle~$\omega_{QR_i}$ in a proper direction.

More generally, we can consider going around edge~$QR_i$ but starting from some arbitrary ``remote'' tetrahedron~$\mathcal T \subset \St(Q)$ and going f\/irst from~$\mathcal T$ to~$P_1P_2R_iQ$ through some two-dimen\-sion\-al faces in the triangulation, then going around~$QR_i$ as in the preceding paragraph, and f\/inally returning from~$P_1P_2R_iQ$ to~$\mathcal T$ following the same way. The corresponding transformation from the old coordinate system to the new one will be again an orthogonal rotation around some axis going through~$Q$.

Now, imagine we are within some chosen tetrahedron~$\mathcal T$ belonging to the star of~$Q$ (i.e., having $Q$ as one of its vertices). For clarity, let us be at a point~$N$ on a small sphere around~$Q$. We are going to draw some special closed paths on this sphere minus $m$~punctures~-- the points where the sphere intersects the edges~$QR_i$. It is clear that we can choose a set of paths $\alpha_i$, $i=1,\dots,m$, in such punctured sphere, with the following properties:
\begin{itemize}
\itemsep=0pt
	\item each $\alpha_i$ begins and ends in~$N$,
	\item within each $\alpha_i$ lies exactly one puncture, namely the intersection of sphere with~$QR_i$, and $\alpha_i$ goes around it in the counterclockwise direction,
	\item the composition $\alpha_1 \circ \dots \circ \alpha_m$ is homotopic to the trivial path in the punctured sphere staying at the point~$N$ .	
\end{itemize}

Denote $X_{QR_i}$ the coordinate system transformation corresponding to~$\alpha_i$. Choosing vertex~$Q$ as the origin of coordinates, we can represent all $X_{QR_1}, \ldots, X_{QR_m}$ as elements of the group~${\rm SO}(3)$. The collapsibility of the composition of~$\alpha_i$ leads to the relation
\begin{equation}
X_{QR_m} \circ \dots \circ X_{QR_1} = \mathop{\rm Id}
\label{eq_rel_in_SO3}
\end{equation}
(the inverse order of $X_{QR_i}$ compared to~$\alpha_i$ ref\/lects the fact that, in the product in~\eqref{eq_rel_in_SO3}, the rightmost elements comes f\/irst, while in the product of~$\alpha_i$~-- traditionally the leftmost).

Now we turn to the inf\/initesimal case. Here, up to second-order inf\/initesimals, $X_{QR_i}$ are just rotations around corresponding edges~$QR_i$ (regardless of how exactly the path~$\alpha_i$ goes between the other edges):
\begin{displaymath}
X_{QR_i} = \mathop{\rm Id} + x_{QR_i}\, d\omega_{QR_i},
\end{displaymath}
where $x_{QR_i}$ is the element of Lie algebra $\mathfrak{so}(3)$ def\/ined by
\begin{displaymath}
(x_{QR_i})_{\alpha\beta} = \sum\limits_\gamma\varepsilon_{\alpha\beta\gamma}\, (\vec e_{QR_i})_\gamma ,
\end{displaymath}
$\varepsilon_{\alpha\beta\gamma}$ being the Levi-Civita symbol.
Equation~(\ref{eq_rel_in_SO3}) gives
\begin{displaymath}
\sum_{i=1}^m x_{QR_i}\, d\omega_{QR_i} = 0 ,
\end{displaymath}
so we get~(\ref{eq_zero_for_dx*}).
\end{proof}

\subsection[Formulas for $\partial \omega_a/\partial l_b$]{Formulas for $\boldsymbol{\partial \omega_a/\partial l_b}$}

The main ingredient of our complex~(\ref{eq_alg_complex}) is matrix $f_3 = \big(\frac{\partial \omega_a}{\partial l_b}\big)$, where $a$ and $b$ run over all edges in the triangulation and the partial derivatives are taken at such values of all lengths~$l_b$ where all $\omega_a = 0$. We are going to express these derivatives in terms of edge lengths and tetrahedron volumes. Recall that all the tetrahedra we deal with are supposed to be consistently oriented (see Convention~\ref{conv2}), so that their volumes have signs. Nonzero derivatives are obtained in the following cases.

\subsubsection*{1st case}
Edges $a=MN$ and $b=QR$ are skew edges of tetrahedron $MNRQ$:
\begin{equation}
\frac{\partial \omega_{MN}}{\partial l_{QR}} = -\frac {l_{MN}\,l_{QR}}{6}\, \frac {1}{V_{MNRQ}}.
\label{eq MNQR}
\end{equation}
This formula is an easy exercise in elementary geometry.

\subsubsection*{2nd case}
Edges $a=MN$ and $b=MP$ belong to the common face of two neighboring tetrahedra $MNPQ$ and $RMNP$ (see the left-hand side of Fig.~\ref{fig2a}):
\begin{equation}
\frac{\partial \omega_{MN}}{\partial l_{MP}} = \frac {l_{MN}\,l_{MP}}{6}\,\frac {V_{NPRQ}}{V_{MNPQ}\,V_{RMNP}}.
\label{eq MNMP}
\end{equation}
To prove this formula, suppose the lengths $l_{MP}$ and $l_{QP}$ are free to change, while the lengths of the remaining seven edges in the two mentioned tetrahedra, \emph{and also the length~$l_{QR}$}, are f\/ixed (this implies of course f\/ixing the dihedral angles at edge~$MN$). Then $l_{QP}$ is a function of~$l_{MP}$ and, using formulas of type~(\ref{eq MNQR}), one can f\/ind that
\begin{equation}
\frac{\partial l_{QP}}{\partial l_{MP}} = \frac{l_{MP}}{l_{QP}}\, \frac{V_{NPRQ}}{V_{RMNP}}.
\label{eq ll1}
\end{equation}
Then,
\begin{equation*}
0 \equiv \frac{d\omega_{MN}}{dl_{MP}} = \frac{\partial \omega_{MN}}{\partial l_{MP}} + \frac{\partial \omega_{MN}}{\partial l_{QP}}\, \frac{\partial l_{QP}}{\partial l_{MP}} = \frac{\partial \omega_{MN}}{\partial l_{MP}} - \frac {l_{MN}\,l_{MP}}{6}\,\frac {V_{NPRQ}}{V_{MNPQ}\,V_{RMNP}}
\end{equation*}
and formula~(\ref{eq MNMP}) follows.

\subsubsection*{3rd case}
Edge $a=b=QR$ is common for exactly three tetrahedra $MNRQ$, $NPRQ$ and~$PMRQ$ (as in the right-hand side of Fig.~\ref{fig2a}):
\begin{equation}
\frac{\partial \omega_{QR}}{\partial l_{QR}} = -\frac {l_{QR}^2}{6}\, \frac {V_{MNPQ}\,V_{RMNP}}{V_{MNRQ}\,V_{NPRQ}\,V_{PMRQ}}.
\label{eq QRQR}
\end{equation}
Observe that this is the most important formula allowing us to construct manifold invariants based on three-dimensional Euclidean geometry. To prove it, consider the ten edges in the three mentioned tetrahedra and suppose that the lengths $l_{MP}$ and $l_{QR}$ are free to change, while the lengths of the remaining eight edges are f\/ixed. Then,
\begin{equation}
\frac{\partial l_{MP}}{\partial l_{QR}} = -\frac{l_{QR}}{l_{MP}}\, \frac{V_{MNPQ}\,V_{RMNP}}{V_{MNRQ}\, V_{NPRQ}}.
\label{eq ll0}
\end{equation}

To prove~\eqref{eq ll0}, we use the following trick: f\/irst, consider the situation where $l_{QP}$ and~$l_{MP}$ are free to change, and the remaining eight lengths are f\/ixed. Then, $\partial l_{QP} / \partial l_{MP}$ is given by formula~\eqref{eq ll1}. Second, consider the situation where $l_{QP}$ and~$l_{QR}$ are free to change, and the remaining eight lengths are f\/ixed. Now a formula of type~\eqref{eq ll1} gives
\begin{equation}
\frac{\partial l_{QP}}{\partial l_{QR}} = \frac{l_{QR}}{l_{QP}}\, \frac{V_{NQMP}}{V_{MNRQ}}.
\label{eq n}
\end{equation}
Third, let now the \emph{seven} edge lengths except $l_{QP}$, $l_{MP}$ and~$l_{QR}$ be f\/ixed; consider $l_{QP}$ as a~function of the other two and then equate $l_{QP}$ to a constant:
\begin{equation*}
l_{QP}(l_{MP},l_{QR}) = \const .
\end{equation*}
This can be viewed as def\/ining $l_{MP}$ as an \emph{implicit function} of~$l_{QR}$, and its derivative, calculated according to the standard formula and using \eqref{eq ll1} and~\eqref{eq n}, is exactly~\eqref{eq ll0}.

Now we can write (assuming again that only $l_{MP}$ and $l_{QR}$ can change):
\begin{equation*}
0 \equiv \frac{d\omega_{QR}}{dl_{QR}} = \frac{\partial \omega_{QR}}{\partial l_{QR}} + \frac{\partial \omega_{QR}}{\partial l_{MP}}\, \frac{\partial l_{MP}}{\partial l_{QR}} = \frac{\partial \omega_{QR}}{\partial l_{QR}} + \frac {l_{QR}^2}{6}\, \frac{V_{MNPQ}\, V_{RMNP}}{V_{MNRQ}\, V_{NPRQ}\, V_{PMRQ}} \, ,
\end{equation*}
and formula~(\ref{eq QRQR}) follows.

\subsubsection*{4th case}
Edge $a=b=QR$ is common for $n > 3$ tetrahedra $P_{i-1}P_iRQ$ $(i = 1,\ldots,n)$:
\begin{equation}
\frac{\partial \omega_{QR}}{\partial l_{QR}} = -\frac {l_{QR}^2}{6}\, \sum\limits_{i=3}^n
\frac {V_{P_1P_{i-1}P_iQ}\,V_{RP_1P_{i-1}P_i}}{V_{P_1P_{i-1}RQ}\,V_{P_{i-1}P_iRQ}\,V_{P_iP_1RQ}}.
\label{eq QRQRn}
\end{equation}
This formula comes out if we draw $n-3$ diagonals $P_1P_i$ ($i = 3,\ldots,n-1$) and apply~(\ref{eq QRQR}) to each of f\/igures $P_1P_{i-1}P_iQR$. Adding up the def\/icit angles around $QR$ in each of those f\/igures and cancelling out the dihedral angles which enter twice and with the opposite signs, we get nothing else than $\omega_{QR}$ for the whole f\/igure.

\section{The acyclic complex and the invariant}
\label{sec_acyclic_general}

\subsection{Generalities on acyclic complexes and their Reidemeister torsions}

We brief\/ly review basic def\/initions from the theory of based chain complexes and their associated Reidemeister torsions, see monograph~\cite{Turaev:2001} for details.

\begin{definition}
\label{dfn:cmplx}
Let $C_0, C_1, \ldots , C_n$ be f\/inite-dimensional $\mathbb R$-vector spaces. The sequence of vector spaces and linear mappings
\begin{equation}
\label{eq:cmplx}
C_* = 0 \xrightarrow{f_n} C_n \xrightarrow{f_{n-1}} C_{n-1} \xrightarrow{} \cdots \xrightarrow{} C_1 \xrightarrow{f_0} C_0 \xrightarrow{f_{-1}} 0
\end{equation}
is called a \emph{chain complex} if $\Ima f_i \subset \Ker f_{i-1}$ for all $i = 0, \ldots, n$. This condition is equivalent to $f_{i-1} \circ f_i = 0$.
\end{definition}

Suppose that chain complex $C_*$ is \emph{based}, i.e., each $C_i$ is endowed with a distinguished basis. To make the notations in this paper consistent with our previous papers, we denote this basis~$\mathcal C_i$ (rather than~$\mathbf c_i$, which is the usual notation in the literature on torsions). The linear mapping $f_i \colon C_{i + 1} \to C_i$ can thus be identif\/ied with its matrix.

\begin{definition}
\label{dfn:acycl}
The based chain complex $C_*$ is said to be \emph{acyclic} if $\Ima f_i = \Ker f_{i-1}$ for all~$i = 0, \ldots, n$.
\end{definition}

\begin{remark}
This condition is equivalent to $\rank f_{i-1} = \dim C_i - \rank f_i$. Since $\Ima f_n = 0$ and $\Ker f_{-1} = C_0$, it follows that in an acyclic complex $f_{n-1}$ is injective and~$f_0$ is surjective.
\end{remark}

Suppose that the chain complex $C_*$ def\/ined in~\eqref{eq:cmplx} is acyclic. For every $i = 0, \ldots, n$, let $\mathcal B_i$ be a subset of the basis~$\mathcal C_i$ such that $f_{i-1} (\mathcal B_i)$ is a basis for~$\Ima f_{i-1}$. Recall that we consider the linear operator~$f_i$ as a matrix whose columns and rows correspond to the basis vectors in~$C_{i+1}$ and~$C_i$ respectively. Denote by $f_i|_{\mathcal B_{i+1}, \overline{\mathcal B}_i}$ the submatrix of~$f_i$ consisting of columns corresponding to $\mathcal B_{i+1}$ and rows corresponding to~$\overline{\mathcal B}_i = \mathcal C_i \setminus \mathcal B_i$.

Due to the acyclicity, $\#\mathcal B_{i+1} = \rank f_i = \dim C_i - \rank f_{i-1} = \dim C_i - \#\mathcal B_i$ and it follows that $f_i|_{\mathcal B_{i+1}, \overline{\mathcal B}_i}$ is square. It is also nondegenerate: its determinant coincides with the determinant of the transition matrix between the two bases $\mathcal C_i$ and $\mathcal B_i \cup f_i(\mathcal B_{i+1})$ for~$C_i$.

\begin{definition}
\label{dfn:tors}
The sign-less quantity
\begin{equation}
\label{eq:torsion}
\tau(C_*) = \prod\limits_{i = 0}^{n-1} \bigl(\det f_i|_{\mathcal B_{i+1}, \overline{\mathcal B}_i}\bigr)^{(- 1)^{i+1}} \in \mathbb{R}^* / \{\pm 1\}
\end{equation}
is called the \emph{Reidemeister torsion} of the acyclic based chain complex~$C_*$.
\end{definition}

\begin{remark}[\cite{Turaev:2001}]
It is easy to show that $\tau(C_*)$ does not depend on the choice of the subsets~$\mathcal B_i$, but of course depends on the choice of the bases in~$C_i$'s.
\end{remark}

\begin{remark}\label{rem1B}
We will also use simplif\/ied notations for~$f_i|_{\mathcal B_{i+1}, \overline{\mathcal B}_i}$, explaining their meaning in the text, as in equation~\eqref{eq_tau} below.
\end{remark}

\subsection{The acyclic complex}

Consider the following sequence of vector spaces and linear mappings:
\begin{equation}
0\xrightarrow{} (dx)'\xrightarrow{f_2} (dl)' \xrightarrow{f_3=f_3^{\rm T}}(d\omega)' \xrightarrow{f_4=-f_2^{\rm T}} \oplus' \mathfrak{so}(3) \xrightarrow{} 0\,.
\label{eq_alg_complex}
\end{equation}
Here is the detailed description of the vector spaces in the chain complex~(\ref{eq_alg_complex}):
\begin{itemize}
\itemsep=0pt
  \item the f\/irst vector space $(dx)'$ is the vector space spanned by the dif\/ferentials of coordinates of all vertices except $A$, $B$, $C$ and~$D$;
  \item the second vector space $(dl)'$ is the vector space spanned by the dif\/ferentials of edge lengths for all edges except those depicted in Figs.~\ref{fig1a} and~\ref{fig1b};
  \item similarly, the third vector space $(d\omega)'$ is the vector space spanned by the dif\/ferentials of def\/icit angles corresponding to the same edges;
  \item the last vector space is a direct sum of copies of the Lie algebra $\mathfrak{so}(3)$ corresponding to the same vertices in the triangulation as in the f\/irst space.
\end{itemize}

Before giving the detailed def\/initions of mappings $f_2$, $f_3$ and~$f_4$, here are some remarks.
\begin{remark}
We use the notations ``$f_2$'' and ``$f_3$'' to make them consistent with other papers on the subject, such as, e.g.~\cite{M,M-thesis}. So, the reader must not be surprised with not f\/inding any~``$f_1$'' in this paper.
\end{remark}
\begin{remark}
\label{rem-b}
There is a natural basis in each of the vector spaces, given by the corresponding dif\/ferentials in $(dx)'$, $(dl)'$ and~$(d\omega)'$, and by the standard Lie algebra generators in~$\oplus' \mathfrak{so}(3)$. Such basis is determined up to an ordering of the vertices in $(dx)'$ and~$\oplus' \mathfrak{so}(3)$, and up to an ordering of the edges in $(dl)'$ and~$(d\omega)'$. Thus, we can identify the elements of vector spaces with column vectors, and mappings~-- with matrices.

For example, the vector space~$(dx)'$ consists of columns of the kind
\begin{displaymath}
\big(dx_{E_1}, dy_{E_1}, dz_{E_1}, \ldots, dx_{E_{N'_0}}, dy_{E_{N'_0}}, dz_{E_{N'_0}}\big)^{\rm T},
\end{displaymath}
where $E_1,\ldots,E_{N'_0}$ are all the vertices in the triangulation except $A$, $B$, $C$ and~$D$. We use thus notation~$N'_0$ for the number of vertices that are \emph{inner} for the manifold with boundary ``$M$~minus the interior of two tetrahedra~$ABCD$ (and with doubled edges $AD$ and~$BC$)'', which is consistent with our other papers, where simply~$N_0$ denotes the number of \emph{all} vertices. Note also that the special role of the edges depicted in Figs.~\ref{fig1a} and~\ref{fig1b}~-- boundary edges for the mentioned manifold, if we double $AD$ and~$BC$~-- consists in the fact that they simply do not take part in forming the second and third linear spaces in complex~(\ref{eq_alg_complex}).
\end{remark}
\begin{remark}
The superscript~$\rm T$ means matrix transposing; the equalities over the arrows in complex~(\ref{eq_alg_complex}) will be proved soon after we def\/ine the mappings $f_2$, $f_3$ and~$f_4$, see Theorem~\ref{th_symmetry}.
\end{remark}
\begin{remark}
The fact that the $f$'s are numbered in the decreasing order in~\eqref{eq:cmplx} and in the increasing order in~\eqref{eq_alg_complex} brings no dif\/f\/iculties.
\end{remark}
\begin{remark}
In our case, there will be no problems with the sign of the torsion, and thus the factoring by~$\{ \pm 1 \}$, as in~\eqref{eq:torsion}, will be unnecessary, see Subsection~\ref{subsec-torsion-invariant} for details.
\end{remark}

Here are the def\/initions of the mappings in the chain complex~(\ref{eq_alg_complex}):
\begin{itemize}
\itemsep=0pt
\item the def\/inition of mapping $f_2$ is obvious: if we change inf\/initesimally the coordinates of vertices, then the corresponding edge length changes are obtained by dif\/ferentiating formulas of the kind
\begin{equation}
l_{MN} = \sqrt{(x_N-x_M)^2+(y_N-y_M)^2+(z_N-z_M)^2},
\label{eq_1a}
\end{equation}
where $M$ and $N$ are two vertices, $x_M,\ldots,z_N$~-- their coordinates, and $l_{MN}$~-- the length of edge~$MN$;
\item mapping $f_3$ goes according to formulas~(\ref{eq MNQR})--(\ref{eq QRQRn});
\item for the mapping~$f_4$, the element of the Lie algebra corresponding to a given vertex, arising from given curvatures~$d\omega$ due to~$f_4$, is by def\/inition given by the left-hand side of formula~(\ref{eq_zero_for_dx*}).
\end{itemize}

\begin{theorem}
Sequence \eqref{eq_alg_complex} is a chain complex, i.e., the composition of two successive mappings is zero.
\label{th2}
\end{theorem}

\begin{proof}
The equality $f_3 \circ f_2 = 0$ is obvious from geometric considerations. Indeed, the edge length changes caused by changes of vertex coordinates give no def\/icit angles, because the whole picture (vertices and edges) does not go out of the Euclidean space~$\mathbb R^3$.

By Lemma~\ref{lemma_infinitesimal}, the equality $f_4 \circ f_3 = 0$ holds as well.
\end{proof}

Complex~(\ref{eq_alg_complex}) can be called a \emph{complex of infinitesimal geometric deformations}. It turns out to have an interesting symmetry property.
\begin{theorem}\label{th_symmetry}
The matrices of mappings in complex \eqref{eq_alg_complex} satisfy the following symmetry properties:
\begin{equation}
f_3 = f_3^{\rm T},\qquad f_4 = - f_2^{\rm T}.
\label{eq_symmetry}
\end{equation}
\end{theorem}

\begin{proof} The symmetry of matrix~$f_3$ means that $\partial \omega_b / \partial l_c = \partial \omega_c / \partial l_b$, and it follows directly from explicit formulas~\eqref{eq MNQR} and \eqref{eq MNMP}.

As for the second equality in~(\ref{eq_symmetry}), it can be proved by a direct writing out of matrix elements, i.e., the relevant partial derivatives. For the mapping~$f_2$, one has to dif\/ferentiate relation~(\ref{eq_1a}), and for~$f_4$~-- use the left-hand side of~(\ref{eq_zero_for_dx*}).
\end{proof}

\begin{remark}
There exists also a dif\/ferent proof of the equality $f_3 = f_3^{\rm T}$, which enables us to look at it perhaps from a dif\/ferent perspective, and based on the \emph{Schl\"afli differential identity} for a Euclidean tetrahedron:
\begin{equation}
\sum_{i=1}^6 l_i \, d\varphi_i = 0
\label{eq_Schlaefli_0}
\end{equation}
for any inf\/initesimal deformations ($l_i$ are edge lengths in the tetrahedron, and~$\varphi_i$ are dihedral angles at edges). It follows from~\eqref{eq_Schlaefli_0} that
\begin{equation}
\sum_a l_a \, d\omega_a = 0,
\label{eq_Schlaefli_1}
\end{equation}
where $a$ runs over all edges in the triangulation.

Consider now the quantity $\Phi=\sum_a l_a \omega_a$ as a function of the lengths~$l_a$, and write the following identity for it:
\begin{equation}
\frac{\partial^2 \Phi}{\partial l_b \partial l_c} = \frac{\partial^2 \Phi}{\partial l_c \partial l_b},
\label{eq_Schlaefli_2}
\end{equation}
where $b$ and $c$ are some edges. It is easy to see that~(\ref{eq_Schlaefli_1}) together with~(\ref{eq_Schlaefli_2}) yield exactly the required symmetry.
\end{remark}

The examples below in Sections~\ref{sec_unknot} and~\ref{sec_lenses} show that there are many enough interesting cases where complex~(\ref{eq_alg_complex}) turns out to be acyclic (see Def\/inition~\ref{dfn:acycl}). However, at this time we cannot make this statement more precise.

\begin{convention}
\label{conv3}
From now on, we assume that we are working with an \emph{acyclic complex}.
\end{convention}

\subsection{The Reidemeister torsion and the invariant}
\label{subsec-torsion-invariant}

As complex~(\ref{eq_alg_complex}) is supposed to be acyclic, we associate to it its \emph{Reidemeister torsion} given by
\begin{equation}
\tau = \frac{\bigl(\det f_2 |_{\overline{\mathcal B}}\bigr) \bigl(\det f_4 |_{\overline{\mathcal B}}\bigr)}{\det f_3 |_{\mathcal B}} = (-1)^{N'_0}\, \frac{\bigl(\det f_2 |_{\overline{\mathcal B}}\bigr)^2}{\det f_3 |_{\mathcal B}}
\label{eq_tau}
\end{equation}
(cf. Def\/inition~\ref{dfn:tors} and Remark~\ref{rem1B} after it). Recall that $N'_0$, according to Remark~\ref{rem-b}, is the number of vertices in the triangulation without $A,B,C$ and~$D$. The letter $\mathcal B$ denotes a maximal subset of edges (and thus basis vectors in $(dl)'$ and~$(d\omega)'$; remember that the edges depicted in Figs.~\ref{fig1a} and~\ref{fig1b} have been already withdrawn) for which the corresponding \emph{ diagonal} minor of~$f_3$ does not vanish. We write this minor as $\det f_3|_{\mathcal B}$, where $f_3|_{\mathcal B}$ is the submatrix of $f_3$ whose \emph{rows and columns} correspond to the edges in~$\mathcal B$. The set $\overline{\mathcal B}$ is the complement of~$\mathcal B$ in the set of all edges except those depicted in Figs.~\ref{fig1a} and~\ref{fig1b}, and $f_2|_{\overline{\mathcal B}}$ (resp. $f_4|_{\overline{\mathcal B}}$) is the submatrix of~$f_2$ (resp.~$f_4$) whose \emph{rows} (resp.\ \emph{columns}) correspond to the edges in~$\overline{\mathcal B}$.

\begin{remark}
\label{remark_sign}
As it is known~\cite{Turaev:2001}, usually the Reidemeister torsion is def\/ined up to a sign, so that special measures must be taken for its ``sign-ref\/ining''. This sign is changed when we change the order of basis vectors in any of the vector spaces. In our case, however, this is not a problem: due to the symmetry proved in Theorem~\ref{th_symmetry}, we can choose our torsion in the form~(\ref{eq_tau}) where the numerator is a square and the denominator is a \emph{diagonal} minor. Both thus do not depend on the order of basis vectors.
\end{remark}

We now def\/ine the value
\begin{equation}
I(M) = \tau\, \frac{\prod_{\rm edges}' l^2}{\prod_{\rm tetrahedra}' (-6V)}\, (6V_{ABCD})^4.
\label{eq_invariant}
\end{equation}
Here $\prod_{\rm edges}' l^2$ means the product of squared lengths for all edges except those depicted in Figs.~\ref{fig1a} and~\ref{fig1b}, or simply \emph{inner} edges; $\prod_{\rm tetrahedra}' (-6V)$ means the product of all tetrahedron volumes multiplied by~$(-6)$ except two distinguished tetrahedra~$ABCD$, and $\tau$ is the Reidemeister torsion of complex (\ref{eq_alg_complex}) given by~(\ref{eq_tau}).

\begin{theorem}
Let $T_1$ and $T_2$ be triangulations of the manifold~$M$ with the same chain of two distinguished tetrahedra~$ABCD$ in them. If complex~\eqref{eq_alg_complex} is acyclic for $T_1$, then it is also acyclic for $T_2$, and the value $I(M)$ given by~\eqref{eq_invariant} is the same for $T_1$ and~$T_2$.
\label{th3}
\end{theorem}

\begin{proof}
By Theorem~\ref{th1}, it is enough to show the invariance of~$I(M)$ under relative Pachner moves.

Suppose we are doing a $2\to 3$ Pachner move: two adjacent tetrahedra $MNPQ$ and~$RMNP$ are replaced by three tetrahedra $MNRQ$, $NPRQ$
and $PMRQ$, see Fig.~\ref{fig2a}. Thus, a new edge~$QR$ appears in the triangulation. We are going to express the new matrix~$f_3$ in terms of the old one. Essentially, we follow~\cite[Section~4]{JNMP1}.

Set $\tilde l_{QR} = l_{QR} - l^{(0)}_{QR}$, where $l^{(0)}_{QR}$ is the solution of $\omega_{QR} = 0$, considered as a function of other edge lengths. Then,
\begin{equation}
d\tilde l_{QR} = dl_{QR} - \sum\limits_{i=1}^{N'_1} \frac{\partial l_{QR}}{\partial l_i}{\Big|}_{\omega_{QR} = 0} \, dl_i,
\label{eq_u}
\end{equation}
where $N'_1$ is the number of inner edges (we reserve the notation~$N_1$ for the number of \emph{all} edges, in analogy with our notations $N'_0$ and~$N_0$ for vertices) in the triangulation before doing the move $2\to 3$. As $\tilde l_{QR}=0$ implies $\omega_{QR}=0$, the dif\/ferential $d\omega_{QR}$ depends only on $d\tilde l_{QR}$ and does not depend on the rest of~$dl_i$. This yields
\begin{equation}
\begin{pmatrix}
d\omega_1 \\
\vdots \\
d\omega_{N'_1} \\
d\omega_{QR}
\end{pmatrix} =
\begin{pmatrix}
{} & {} & {} & \frac{\partial\omega_1}{\partial l_{QR}} \\
{} & f_3^{\rm old} & {} & \vdots \\
{} & {} & {} & \frac{\partial\omega_{N'_1}}{\partial l_{QR}} \\
0 & \cdots & 0 & \frac{\partial\omega_{QR}}{\partial l_{QR}}
\end{pmatrix}
\begin{pmatrix}
dl_1 \\
\vdots \\
dl_{N'_1} \\
d\tilde l_{QR}
\end{pmatrix},
\label{eq_v}
\end{equation}
where $f_3^{\rm old}$ is, of course, the matrix~$f_3$ before the move. Formulas \eqref{eq_v} and~\eqref{eq_u} give the representation of the ``new'' matrix $f_3^{\rm new}$ as the following product:
\begin{equation*}
f_3^{\rm new} =
\begin{pmatrix}
{} & {} & {} & \frac{\partial\omega_1}{\partial l_{QR}} \\
{} & f_3^{\rm old} & {} & \vdots \\
{} & {} & {} & \frac{\partial\omega_{N'_1}}{\partial l_{QR}} \\
0 & \cdots & 0 & \frac{\partial\omega_{QR}}{\partial l_{QR}}
\end{pmatrix}
\begin{pmatrix}
1 & {} & {} & 0 \\
{} & \ddots & {} & \vdots \\
{} & {} & 1 & 0 \\
- \frac{\partial l_{QR}}{\partial l_1} & \cdots & - \frac{\partial l_{QR}}{\partial l_{N'_1}} & 1
\end{pmatrix}.
\end{equation*}

The ``new'' set of edges~$\mathcal B^{\rm new}$ in~\eqref{eq_tau} can be taken as $\mathcal B^{\rm new} = \mathcal B^{\rm old} \cup\{QR\}$, then the minors of $f_2$ and~$f_4$ remain the same. Using~\eqref{eq QRQR}, we obtain the ratio between the new and old minors of~$f_3$:
\begin{equation}
\frac{\det f_3^{\rm new} |_{\mathcal B^{\rm new}}}{\det f_3^{\rm old} |_{\mathcal B^{\rm old}}} = \frac{\partial\omega_{QR}}{\partial l_{QR}} = - \frac{l_{QR}^2}{6}\frac{V_{MNPQ}V_{RMNP}}{V_{MNRQ}V_{NPRQ}V_{PMRQ}}  .
\label{eq QRQR minors}
\end{equation}
This together with \eqref{eq_invariant} and~\eqref{eq_tau} proves that $I(M)$ does not change under $2\to 3$ and $3\to 2$ Pachner moves.

Now we consider a $1\to 4$ Pachner move. It means that a tetrahedron $MNPQ$ is divided into four tetrahedra by adding a new vertex~$R$ inside it, as in Fig.~\ref{fig2b}. Hence, three new components are added to the vectors in the space~$(dx)'$~-- the f\/irst vector space in sequence~(\ref{eq_alg_complex}), namely the dif\/ferentials $dx_R$, $dy_R$ and~$dz_R$ of coordinates of vertex~$R$. In the same way, four new components are added to the vectors in the space~$(dl)$~-- the second vector space in~(\ref{eq_alg_complex}), namely $dl_{MR}$, $dl_{NR}$, $dl_{PR}$ and~$dl_{QR}$. We add the edge $QR$ to the set~$\mathcal B$, then $MR$, $NR$ and~$PR$ are added to~$\overline{\mathcal B}$. The minor of $f_2$ acquires a block triangular structure and becomes the product of the old minor by the determinant of a new $3\times 3$ block, namely
\begin{equation*}
\frac{dl_{MR}\wedge dl_{NR}\wedge dl_{PR}}{dx_R\wedge dy_R\wedge dz_R} = \frac{6V_{MNPR}}{l_{MR}\,l_{NR}\,l_{PR}}
\end{equation*}
(this equality follows from elementary geometry, compare formulas (31) and~(32) in~\cite{JNMP1}). Due to the same considerations as in the case of a $2\to 3$ Pachner move, the minor of $f_3$ gets multiplied by the very same factor given in~(\ref{eq QRQR}). Comparing that with equations~(\ref{eq_tau}) and~(\ref{eq_invariant}), we see that our value~$I(M)$ does not change under a $1\to 4$ Pachner move, as well as under the inverse move.

The fact that no minors considered in our proof vanished obviously implies that the acyclicity of our complex~(\ref{eq_alg_complex}) is preserved under the Pachner moves.
\end{proof}

Recall that in the beginning of Subsection~\ref{subs_daav}, we introduced a function (let us denote it now by $\phi$) ascribing to the vertices of a triangulation of manifold~$M$ three real numbers in arbitrary way with the only condition: the volumes of all tetrahedra calculated by~(\ref{eq_oriented_volume}) must be nonzero. Let us call $\phi$ a \emph{geometrization function} on the set of vertices. To ensure the ``well-def\/inedness'' of our invariant, we must show that it is independent of the choice of the geometrization function.
\begin{theorem}
Let $\phi_1$ and $\phi_2$ be geometrization functions such that $\phi_1|_{\{A, B, C, D\}} = \phi_2|_{\{A, B, C, D\}}$, where $A$, $B$, $C$ and~$D$ belong to the chain of two distinguished tetrahedra. If $I_{\phi_1}(M)$ and $I_{\phi_2}(M)$ are the corresponding values of the invariant~\eqref{eq_invariant}, then $I_{\phi_1}(M) = I_{\phi_2}(M)$.
\end{theorem}

\begin{proof}
Let $P$ be any vertex in the triangulation, not belonging to the chain of two distinguished tetrahedra. First, we are going to prove that there exists a sequence of relative Pachner moves removing~$P$ from the triangulation.

We begin with subdividing the initial triangulation so that it becomes combinatorial. Denote by $\St(P)$ and $\lk(P)$ the star and the link of~$P$ respectively in this combinatorial triangulation.

Let two adjacent tetrahedra $MNPQ$ and $RMNP$ belong to~$\St(P)$. First, we do a move $2 \to 3$ by adding a new edge~$QR$. Then, the tetrahedron $MNRQ$ goes out of $\St(P)$ and the edge $MN$ goes out of~$\lk(P)$. Continuing this way for other pairs of adjacent tetrahedra from~$\St(P)$ containing $M$ as a vertex, we can reduce to three the number of edges starting from $M$ and belonging to~$\lk(P)$. After that, doing a move $3 \to 2$, we remove the edge $MP$ from $\St(P)$ and the vertex $M$ from~$\lk(P)$. Continuing this way, we reduce to four the number of vertices in~$\lk(P)$.

Now, we do a move $4 \to 1$ to remove $P$ from the triangulation. After that, doing a move $1 \to 4$, we return it back, but with any other ascribed coordinates. Finally, we return to the initial triangulation inverting the whole sequence of Pachner moves. It remains to say that, due to the invariance under Pachner moves, the value $I(M)$ is unchanged at every step.
\end{proof}

\begin{remark}
Invariant $I(M)$ in~(\ref{eq_invariant}) can depend a priori on the geometry of the distinguished tetrahedron~$ABCD$. However, it turns out that with the multiplier $(6V_{ABCD})^4$ introduced in~(\ref{eq_invariant}) the invariant is just a number at least in the examples considered below in
Sections~\ref{sec_unknot} and~\ref{sec_lenses}.
\end{remark}

\section{How to change the framing}
\label{sec_framing}

Just as Pachner moves are elementary rebuildings of a triangulation of a \emph{closed} manifold, \emph{shellings} and inverse shellings are elementary rebuildings of a triangulation of a manifold \emph{with boundary}, see~\cite[Section~5]{L}. A topological f\/ield theory dealing with triangulated manifolds must answer the question what happens with an invariant like~$I(M)$ under shellings.

While we leave a general answer to this questions to further papers, we will explain in this section how some shellings on the toric boundary of our manifold ``$M$~minus two tetrahedra'' correspond to changing the framing of the knot determined by these two tetrahedra. We also show what happens with matrices $f_3$ and~$f_2$ from complex~(\ref{eq_alg_complex}) under these shellings. These results will be used in calculations in Sections~\ref{sec_unknot} and~\ref{sec_lenses}.

\begin{figure}[ht]
\centerline{\includegraphics[scale=0.4]{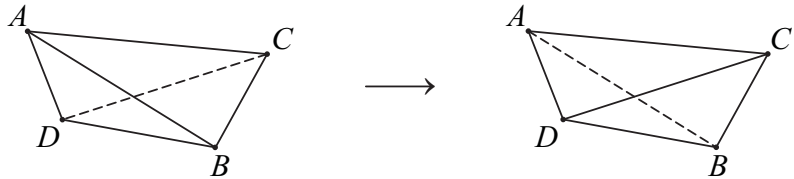}}
\caption{Turning a tetrahedron inside out, thus changing the framing by $1/2$.}
\label{fig4}
\end{figure}

It is enough to show how to change the knot framing by one-half of a revolution. We can achieve this if we manage to ``turn inside out'' one of the tetrahedra~$ABCD$ in Figs.~\ref{fig1a} and~\ref{fig1b}, e.g., in the way shown in Fig.~\ref{fig4}.
\begin{remark}
Of course, the framing can be changed in other direction similarly. In this case, we should f\/irst draw the left-hand-side tetrahedron in Fig.~\ref{fig4} as viewed from another direction, so as the diagonals of its projection are $AC$ and $BD$, instead of $AB$ and $CD$ in Fig.~\ref{fig4}. Then we replace the dashed ``diagonal'' with the solid one and vice versa.
\end{remark}

Return to Fig.~\ref{fig4}. In order to be able to glue the ``turned inside out'' tetrahedron back into the triangulation, we can glue to it two more tetrahedra~$ABCD$: one to the front and one to the back. So, we glue the same tetrahedron as drawn in the left-hand side of Fig.~\ref{fig4}, to the \emph{two} ``front'' faces, $ADC$ and $DBC$, of the ``turned inside out'' tetrahedron in the right-hand side of Fig.~\ref{fig4}, and again the same tetrahedron as in the left-hand side of Fig.~\ref{fig4} to the two ``back'' faces, $ABC$ and $ADB$ (always gluing a vertex to the vertex of the same name). After this, the obtained ``sandwich'' of three tetrahedra can obviously be glued into the same place which was occupied by the single tetrahedron in the left-hand side of Fig.~\ref{fig4}.

How will the invariant~$I(M)$ change? Of course, the product of tetrahedron volumes in~(\ref{eq_invariant}) will be multiplied by $(6V_{ABCD})^2$, and the product of edge lengths will be multiplied by~$l_{AB}$ and~$l_{CD}$, because we have added two tetrahedra $ABCD$ to our triangulation, and the edges~$AB$ and~$CD$ of the initial tetrahedron in the left-hand side of Fig.~\ref{fig4} changed their status from being inner to lying on the boundary.

To describe the change of matrix $f_3$, it is f\/irst convenient to introduce matrix~$F_3$, which consists by def\/inition of \emph{all} partial derivatives $(\partial\omega_a / \partial l_b)$, \emph{including the edges that belong to the distinguished tetrahedra} in Figs.~\ref{fig1a} and~\ref{fig1b}. Thus, $f_3$ is a submatrix of~$F_3$. Then we introduce a~``normalized'' version of matrix~$F_3$, denoted~$G_3$, as follows:
\begin{equation}
G_3 = 6V_{ABCD}\, \diag\big(l_1^{-1},\ldots,l_{N_1}^{-1}\big)\,F_3\,\diag\big(l_1^{-1},\ldots,l_{N_1}^{-1}\big).
\label{eq_*}
\end{equation}
Here $N_1$ is the total number of edges in the triangulation of the manifold~$M$. Just as $f_3$, mat\-ri\-ces~$F_3$ and~$G_3$ are symmetric.

Now we describe what happens with $G_3$ when we change the framing. We represent $G_3$ in a~block form where the last row and the
last column correspond to the edge~$CD$, and the next to last row and column to the edge~$AB$:
\begin{equation}
G_3=\left( \begin{array}{cc} K & L^{\rm T} \\[.5ex] L & \begin{array}{cc} \alpha & \beta \\
\beta & \gamma \end{array} \end{array}\right). \label{eq_**}
\end{equation}
Here $\alpha$, $\beta$ and~$\gamma$ are real numbers, $K$ is an $(N_1-2)\times (N_1-2)$ block, and $L$ is a $2\times (N_1-2)$ block.

Recall that we have chosen a consistent orientation for all tetrahedra in the triangulation, which means, for every tetrahedron, a proper ordering of its vertices up to even permutations. The initial tetrahedron in the left-hand side of Fig.~\ref{fig4} thus can have either orientation $ABCD$ or~$BACD$. Let $\epsilon=1$ for the f\/irst case and $\epsilon=-1$ for the second case.

\begin{theorem}
\label{th_sec4}
Let $\widetilde G_3$ be the matrix $G_3$ after the change of framing which adds two new edges $\widetilde{AB}$ and $\widetilde{CD}$ to the triangulation in the way described above. Then,
\begin{equation}
\widetilde G_3 = \left( \begin{array}{ccc} K & L^{\rm T} & {\bf 0} \\[.5ex]
L & \begin{array}{cc} \alpha & \beta-\epsilon \\ \beta-\epsilon & \gamma \end{array} & \begin{array}{cc} 0 & \;\,\epsilon\;\, \\
\;\,\epsilon\;\, & 0 \end{array} \\ {\bf 0} & \begin{array}{cc} 0 & \;\;\;\epsilon\;\;\; \\ \;\;\;\epsilon\;\;\; & 0
\end{array} & \begin{array}{cc} 0 & -\epsilon \\ -\epsilon & 0 \end{array}\end{array} \right),
\label{eq_***}
\end{equation}
where the two last rows and columns correspond to the new added edges $\widetilde{AB}$ and~$\widetilde{CD}$.
\end{theorem}

\begin{proof}
The normalization (\ref{eq_*}) of matrix $G_3$ has been chosen keeping in mind formulas of type~(\ref{eq MNQR}). The derivatives like $\partial \omega_{AB}/\partial l_{CD}$ contribute to the elements of~$G_3$ as $-1$ if the orientation of the corresponding tetrahedron is~$ABCD$, and as~$+1$ otherwise. This, f\/irst, explains why~$\epsilon$ is subtracted from the matrix element~$\beta$ when the initial edges $AB$ and~$CD$ cease to belong to the same tetrahedron. Second, it explains the appearance of $\pm\epsilon$ in the last two rows and columns of~$\widetilde G_3$. It remains to show that the other new matrix elements vanish, e.g.,
\begin{equation}
\frac{\partial \omega_{AC}}{\partial l_{\widetilde{AB}}} = 0, \qquad \frac{\partial \omega_{\widetilde{AB}}}{\partial l_{\widetilde{AB}}} = 0 \quad
\hbox{and so on,}
\label{eq_*****}
\end{equation}
and that some other elements in $G_3$ do not change while, seemingly, the triangulation change has touched them.

The f\/irst equality in~(\ref{eq_*****}) is due to the fact that $l_{\widetilde{AB}}$ inf\/luences two dihedral angles which enter in $\omega_{AC}$;
these angles belong to two tetrahedra $ABCD$ which dif\/fer only in their orientations and thus sum up to an identical zero. A similar explanation
works for the second equality in~(\ref{eq_*****}) as well. Moreover, a similar reasoning shows that, although new summands are added to some elements in $G_3$ like $\partial \omega_{AC} / \partial l_{AD}$, these summands cancel each other because they belong to tetrahedra with opposite
orientations.
\end{proof}

From matrix $\widetilde G_3$, we can obtain the new matrix $F_3$, and then take its relevant submatrix as new~$f_3$. As for the matrix~$f_2$
in complex~(\ref{eq_alg_complex}), it will just acquire two new rows whose elements, like all elements in~$f_2$, are obtained by dif\/ferentiating relations of type~(\ref{eq_1a}).

\section{Calculation for unknot in three-sphere}
\label{sec_unknot}

We now turn to applications of our ideas to examples~-- manifolds with framed knots in them. In this section we do calculations for the simplest case~-- an unknot in a three-sphere.

Let $K$ be a knot in three-sphere, then we denote $N(K)$ its open tubular neighbourhood and $M_K = S^3 \setminus N(K)$ its exterior. Recall that $M_K$ is a three-manifold with the boundary consisting of a single torus.

\begin{figure}[ht]
\centerline{\includegraphics[scale=0.3]{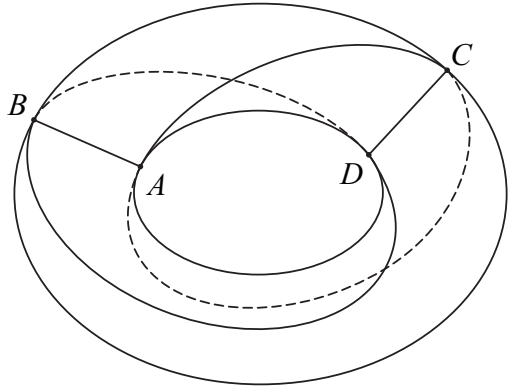}}

\caption{Beginning of the construction of the triangulated solid torus: a chain of two tetrahedra $ABCD$.}
\label{fig1a-modified}
\end{figure}

From now on we suppose that $K = \bigcirc$ is unknot. Then, $M_\bigcirc$ is a f\/illed torus. We glue $M_\bigcirc$ out of six tetrahedra in the following way. First, we take two identically oriented tetrahedra $ABCD$ and glue them together in much the same way as in Fig.~\ref{fig1a}, but using edges $AB$ and~$CD$ for gluing, see Fig.~\ref{fig1a-modified}.
What prevents this chain from being a f\/illed torus is its ``zero thickness'' at the edges $AB$ and~$CD$. We are going to remove this dif\/f\/iculty by gluing some more tetrahedra to the chain.

First, we glue at ``our'' side one more tetrahedron~$DABC$ (of the opposite orientation!) to faces $ABC$ and~$ABD$. This already creates a nonzero thickness at the edge~$AB$, yet we glue still one more tetrahedron~$ABCD$ to the two free faces of the new tetrahedron~$DABC$, that is, $BCD$ and~$ACD$. In the very same way, in order to remove the zero thickness along edge~$CD$, we glue one more tetrahedron of the opposite orientation at ``our'' side of the f\/igure to faces $BCD$ and~$ACD$, and then glue a tetrahedron $ABCD$ to the two free faces of the new tetrahedron.

To distinguish edges of the same name, we introduce the following notations. Edges $AB$ and~$CD$ present in Fig.~\ref{fig1a-modified} will be denoted $(AB)_1$ and~$(CD)_1$. Then, we think of one of the tetrahedra in Fig.~\ref{fig1a-modified} as f\/irst, and the other as second, and accordingly assign to the rest of their edges indices $1$ or~$2$. It remains to denote four edges, of which two lie inside the f\/illed torus, and two~-- in the boundary. We denote the inner edges as $(AB)_3$ and~$(CD)_3$, and the boundary ones~-- as $(AB)_2$ and~$(CD)_2$.

\begin{figure}[ht]
\centerline{\includegraphics[scale=0.3]{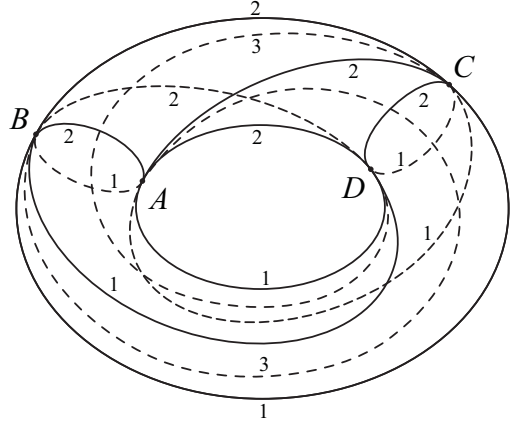}}
\caption{Triangulated f\/illed torus.}
\label{fig8}
\end{figure}

The obtained triangulated f\/illed torus is depicted in Fig.~\ref{fig8}, where the numbers correspond to the subscripts at edges.

Now we can glue the tetrahedron chain from Fig.~\ref{fig1a} to the f\/illed torus in Fig.~\ref{fig8}, gluing together faces of the same name in a natural manner. Thus, sphere~$S^3$ appears with a framed unknot in it determined by the tetrahedra from Fig.~\ref{fig1a}, for which we can calculate the invariant according to Section~\ref{sec_acyclic_general} and then as well the values of the invariant for other framings according to Section~\ref{sec_framing}. We formulate the result as the following
\begin{theorem}
\label{th_Iunknot}
The invariant $I^{(r)}\bigl(M_\bigcirc\bigr)$ for a framing~$r=m$ or $r=m+1/2$, where $m\in \mathbb Z$, is given by{\samepage
\begin{equation*}
I^{(m)}\bigl(M_\bigcirc\bigr) = \frac{1}{m^4}
\qquad \mbox{or}
\qquad
I^{(m + 1/2)}\bigl(M_\bigcirc\bigr) = \frac{- 1}{m^2 \, (m+1)^2},
\end{equation*}
respectively.}
\end{theorem}

\begin{proof}
Our triangulation of sphere $S^3$ does not contain any vertices besides $A$, $B$, $C$ and~$D$. Thus, the spaces $(dx)'$ and $\oplus' \mathfrak{so}(3)$ in complex~\eqref{eq_alg_complex} are zero-dimen\-sional, i.e., complex~\eqref{eq_alg_complex} is reduced to \begin{equation}
0\xrightarrow{} (dl)' \xrightarrow{f_3} (d\omega)' \xrightarrow{} 0,
\label{eq_alg_compl-copy}
\end{equation}
where the dash means that the corresponding quantities are taken only for edges not belonging to the tetrahedra in Fig.~\ref{fig1a}. Formula~\eqref{eq_tau} for the torsion takes thus a simple form
\begin{equation}
\tau = \frac{1}{\det f_3}.
\label{eq_hm_*}
\end{equation}

As we have explained in Section~\ref{sec_framing}, it makes sense to consider matrix~$F_3$ which consists, by def\/inition, of the partial derivatives of \emph{all} def\/icit angles with respect to \emph{all} edge lengths in the triangulation of $S^3$ and of which $f_3$ is a submatrix. Moreover, it makes sense to consider the ``normalized'' version of~$F_3$, i.e., matrix~$G_3$ def\/ined by~(\ref{eq_*}). The invariant~\eqref{eq_invariant} with torsion~\eqref{eq_hm_*} takes then the form
\begin{equation}
I(M_{\bigcirc}) = \frac{(6V_{ABCD})^4 \prod_{\rm edges}' l^2}{\det f_3 \prod_{\rm tetrahedra}' (-6V)} = \frac{1}{\det g_3},
\label{eq_inv0-copy}
\end{equation}
where $g_3$ is the submatrix of $G_3$ consisting of the same rows and columns of which consists $f_3$ as a submatrix of~$F_3$.

Matrix $G_3$ has a simple block structure caused by the vanishing of a derivative $\partial \omega_a / \partial l_b$ in a~case where only pairs of oppositely oriented tetrahedra make contributions in it. Nonzero matrix elements can be present in one of the six blocks corresponding to the following possibilities:
\begin{itemize}
\itemsep=0pt
\item[(i)] $a$ is one of edges $AB$ and $b$ is one of edges $CD$, denote the block of such elements as $S_1$;
\item[(ii)] vice versa: $a$ is one of edges $CD$ and $b$ is one of edges $AB$, these elements form block $S_1^{\mathrm T}$;
\item[(iii)] $a$ is one of edges $AC$ and $b$ is one of edges $BD$, block $S_2$;
\item[(iv)] vice versa, block $S_2^{\mathrm T}$;
\item[(v)] $a$ is one of edges $AD$ and $b$ is one of edges $BC$, block $S_3$;
\item[(vi)] vice versa, block $S_3^{\mathrm T}$.
\end{itemize}
Matrix~$G_3$ has thus the following block structure:
\begin{equation}
G_3 = \left(\begin{array}{c c|c c|c c}{\bf 0}_3 & S_1 &  &  &  &  \\ S_1^{\rm T} & {\bf 0}_3 &  &  &  &  \\\hline  &  & {\bf 0}_2 & S_2 &  &  \\  &  & S_2^{\rm T} & {\bf 0}_2 &  &  \\\hline  &  &  &  & {\bf 0}_2 & S_3 \\  &  &  &  & S_3^{\rm T} &{\bf 0}_2 \end{array}\right).\label{eq_block_structure-copy}
\end{equation}

It is also important that the block structure~\eqref{eq_block_structure-copy} is preserved under the change of~$G_3$ corresponding to a change of framing, as can be checked using~\eqref{eq_**} and~\eqref{eq_***}; although, of course, the sizes of blocks $S_1$ and $S_2$ and consequently their transposes and corresponding zero blocks do change; block~$S_3$ remains intact.

A simple calculation using the explicit form of matrix blocks concludes the proof of the theorem. We think there is no need to present here all details of this calculation, especially because in Section~\ref{sec_lenses} we give details for a similar calculation in a more complicated case of unknots in lens spaces.
\end{proof}

\section{Calculation for unknots in lens spaces}
\label{sec_lenses}

\subsection{Generalities on lens spaces and their triangulations}
\label{subsec_lens_general}

Let $p, q$ be two coprime integers such that $0< q < p$. We identify $S^3$ with the subset $\{(z_1, z_2) \in \mathbb C^2 \colon \; |z_1|^2 + |z_2|^2 = 1\}$ of $\mathbb C^2$. The lens space $L(p,q)$ is def\/ined as the quotient manifold $S^3 / \sim$, where $\sim$ denotes the action of the cyclic group $\mathbb Z_p$ on $S^3$ given by
\begin{equation*}
\zeta \cdot (z_1, z_2) = (\zeta z_1, \zeta^q z_2), \qquad \zeta = e^{2\pi i/p}.
\end{equation*}

As a consequence the universal cover of lens spaces is the three-dimensional sphere $S^3$ and
\begin{equation}
\pi_1 \bigl( L(p,q) \bigr) = H_1 \bigl( L(p,q) \bigr) = \mathbb Z_p.
\label{eq_Z_p}
\end{equation}

The full classif\/ication of lens spaces is due to Reidemeister and is given in the following
\begin{theorem}[\cite{Reidemeister}]
Lens spaces $L(p, q)$ and $L(p', q')$ are homeomorphic if and only if $p' = p$ and $q' = \pm q^{\pm 1} \bmod p$.
\end{theorem}

\begin{figure}[ht]
\centerline{\includegraphics[scale=0.3]{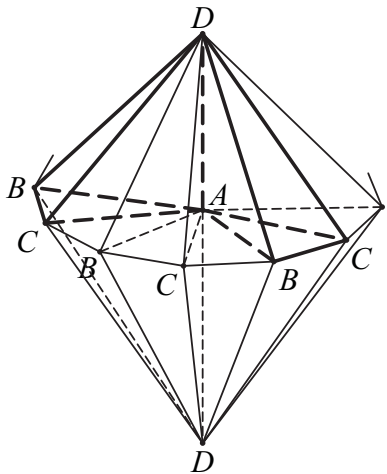}}
\caption{A chain of two tetrahedra in a lens space.}
\label{fig5}
\end{figure}

Now we describe triangulations of $L(p,q)$ which will be used in our calculations. Consider the bipyramid of Fig.~\ref{fig5},
which contains $p$ vertices~$B$ and $p$ vertices~$C$. The lens space $L(p,q)$ is obtained by gluing the upper half of its surface to the lower half, the latter having been rotated around the vertical axis through the angle $2\pi q/p$ in a chosen ``positive'' direction in such way that every ``upper'' triangle $BCD$ is glued to some ``lower'' triangle $BCD$ (the vertices of the same names are identif\/ied).

The connection of Fig.~\ref{fig5} with the above description of $L(p,q)$ as~$S^3 / \sim $ can be established as follows. The bipyramid is identif\/ied with the part of~$S^3$ cut out by the inequalities $0\le \arg z_1 \le 2\pi i/p$, so that $\arg z_1=0$ for points $(z_1,z_2)$ in the upper half of the bipyramid surface and $\arg z_1 = 2\pi i/p$ in its lower half; of course, $z_1=0$ in the equator. Each of these two halves can be parameterized by the complex variable~$z_2$, with $z_2=0$ in the points~$D$ and $|z_2|=1$ in the equator. The rotation of the equator between two nearest points~$B$ in the positive direction corresponds to multiplying $z_2$ by~$\zeta$.

A generator of the fundamental group can be represented, e.g., by some broken line $BCB$ (the two end points~$B$ are dif\/ferent) lying in
the equator of the bipyramid. We assume that a~generator chosen in such way corresponds to the element $1\in \mathbb Z_p$ under the identif\/ication~(\ref{eq_Z_p}).

The boldface lines (solid and dashed) in Fig.~\ref{fig5} single out two \emph{identically oriented} tetrahedra~$ABCD$ which form a chain exactly
like the one in Fig.~\ref{fig1a}. Note that going along this chain (e.g., along the way $BAB$) corresponds to a certain nonzero element from $H_1 \bigl( L(p,q)\bigr) = \mathbb Z_p$. In this paper, we call a knot in $L(p,q)$ determined by a tetrahedron chain of the kind of Fig.~\ref{fig5} an ``unknot'' in~$L(p,q)$. We are going to calculate our invariant for such unknots with dif\/ferent framings.

\subsection[The structure of matrix $(\partial \omega_a/ \partial l_b)$]{The structure of matrix $\boldsymbol{(\partial \omega_a/ \partial l_b)}$}
\label{subsec_lens_structure}

The triangulation of a lens space, described in the previous subsection, does not contain any vertices besides $A$, $B$, $C$ and~$D$. It follows then that complex~(\ref{eq_alg_complex}), corresponding to such triangulation, is reduced to a single mapping~$f_3$, that is, it takes the already known to us form~\eqref{eq_alg_compl-copy}:
\begin{equation}
0\xrightarrow{} (dl)' \xrightarrow{f_3} (d\omega)' \xrightarrow{} 0.
\label{eq_alg_compl}
\end{equation}
This complex is acyclic provided $\det f_3 \ne 0$.

As we have explained in Section~\ref{sec_framing}, it makes sense to consider matrix~$F_3$ which consists, by def\/inition, of the partial
derivatives of \emph{all} def\/icit angles with respect to \emph{all} edge lengths in the triangulation of the lens space and of which $f_3$ is
a submatrix. Moreover, it makes sense to consider the ``normalized'' version of $F_3$, i.e., matrix~$G_3$ def\/ined by~(\ref{eq_*}).

Matrix $G_3$ has many zero entries. They appear in one of two ways: either the corresponding derivative $\partial \omega_a/
\partial l_b$ vanishes because the edges $a$ and~$b$ do not belong to the same tetrahedron, or the cause is like that explained in the proof of Theorem~\ref{th_sec4}, compare~(\ref{eq_*****}).

\begin{figure}[ht]
\centerline{\includegraphics[scale=0.3]{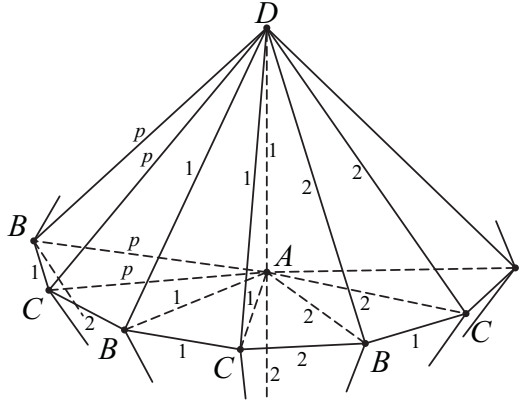}}
\caption{To the explanation of the structure of matrix $G_3$.}
\label{fig6}
\end{figure}

We will denote the triangulation edges by indicating the origin and end vertices of a given edge. As dif\/ferent edges may have the same origin and end vertices, we will assign indices to them, as indicated in Fig.~\ref{fig6}.
For example, as one can see from this f\/igure, there exist $p$ dif\/ferent edges $AB$, and each of them is equipped with an index from 1 to~$p$. So, we denote by $(AB)_n$ the edge $AB$ equipped with index~$n$.

To describe the structure of matrix $G_3$, we introduce the following ordering on the set of all edges in triangulation:
\begin{displaymath}
(AB)_1, \ldots, (AB)_p,\; (CD)_1, \ldots, (CD)_p,
\end{displaymath}
\begin{displaymath}
(AC)_1, \ldots, (AC)_p, \; (BD)_1, \ldots, (BD)_p,
\end{displaymath}
\begin{displaymath}
(AD)_1, (AD)_2,\; (BC)_1, (BC)_2.
\end{displaymath}

In this way, we put in order the basis vectors in spaces $(dl)$ and~$(d\omega)$. The order of matrix $G_3$ is $(4p+4)\times(4p+4)$, and with respect to the preceding ordered basis, $G_3$ has the following block structure:
\begin{equation}
G_3 = \left(\begin{array}{c c|c c|c c}{\bf 0}_p & S_1 &  &  &  &  \\ S_1^{\rm T} & {\bf 0}_p &  &  &  &  \\\hline  &  & {\bf 0}_p & S_2 &  &  \\  &  & S_2^{\rm T} & {\bf 0}_p &  &  \\\hline  &  &  &  & {\bf 0}_2 & S_3 \\  &  &  &  & S_3^{\rm T} &{\bf 0}_2 \end{array}\right),\label{eq_block_structure}
\end{equation}
where $S_1$, $S_2$ are $p \times p$ submatrices and $S_3$ is a $2 \times 2$ submatrix. Here and below the empty spaces in matrices are of course occupied by zeroes.

We now describe the structure of $S_1$, $S_2$ and $S_3$.
\begin{itemize}
\itemsep=0pt
  \item[(i)] The $i$th row of $S_1$ consists of the partial derivatives $\partial \omega_{(AB)_i} / \partial l_{(CD)_j}$, and
with the help of Fig.~\ref{fig6}, we may conclude that there exist exactly four nonzero entries in each row, namely:
\begin{equation}
c\cdot \frac{\partial \omega_{(AB)_i}}{\partial l_{(CD)_i}}, \qquad c\cdot \frac{\partial \omega_{(AB)_i}}{\partial l_{(CD)_{i-1}}}, \qquad c\cdot \frac{\partial
\omega_{(AB)_i}}{\partial l_{(CD)_{i-q}}}, \qquad c\cdot \frac{\partial \omega_{(AB)_i}}{\partial l_{(CD)_{i-q-1}}},
\label{eq_ccdot}
\end{equation}
where
\begin{displaymath}
c=\frac{6V_{ABCD}}{l_{AB}l_{CD}}.
\end{displaymath}
Here, all indices change cyclicly from 1 to $p$, i.e., for instance, $0 \equiv p$, $-1 \equiv p-1$, and so forth. It is convenient to choose the orientation of the four tetrahedra in Figs.~\ref{fig1a} and~\ref{fig1b} as $BACD$. Then, according to~(\ref{eq MNQR}), the expressions in~(\ref{eq_ccdot}) turn respectively into
\begin{displaymath}
1,\quad -1,\quad -1,\quad 1.
\end{displaymath}
Moving further along these lines, we obtain
\begin{equation}
S_1 = {\bf 1}_p - E - E^q + E^{q+1} = ({\bf 1}_p - E^q) ({\bf 1}_p - E),
\label{eq_S1}
\end{equation}
where ${\bf 1}_p$ is the identity matrix of size $p \times p$, and
\begin{equation}
E =\left(\begin{array}{ccccc}0 &  & \ldots &  & 1 \\1 & 0 &  &  &  \\ & 1 & \ddots &  & \vdots \\ &  & \ddots & 0 &  \\ &  &  & 1 & 0\end{array}\right).
\label{eq_E}
\end{equation}
  \item[(ii)] Similarly,
\begin{equation}
S_2 = ({\bf 1}_p - E^q) ({\bf 1}_p - E^{- 1}). \label{eq_S2}
\end{equation}
  \item[(iii)] Finally, one can verify that
\begin{equation*}
S_3 = p \left(
\begin{array}{cc}
1 & - 1 \\
- 1 & 1
\end{array} \right).
\end{equation*}
\end{itemize}

\subsection[Invariant for the ''simplest'' framing]{Invariant for the ``simplest'' framing}
\label{subsec_lens_calculation}

Let $L_n(p,q)$ denote the lens space $L(p,q)$ with a tetrahedron chain like in Fig.~\ref{fig5}, but with the angular distance $n\cdot\frac{2\pi}{p}$ between the two distinguished tetrahedra. Thus, we have an ``unknot going along the element $n\in\mathbb Z_p = H_1 \bigl(L(p,q)\bigr)$'', which we will also call the $n$th unknot. In this section we consider the simplest case when a framed knot is determined directly by a tetrahedron chain of the type depicted in Fig.~\ref{fig5}.

According to~(\ref{eq_tau}) and~(\ref{eq_invariant}) and the form of complex~(\ref{eq_alg_compl}), the invariant takes the simple form
\begin{equation}
I\bigl(L_n(p,q)\bigr) = \frac{1}{\det g_3},
\label{eq_inv0}
\end{equation}
where $g_3$ is the submatrix of $G_3$ consisting of the same rows and columns of which consists $f_3$ as a submatrix of~$F_3$, compare equation~\eqref{eq_inv0-copy}.

We also identify $n\in \mathbb Z_p$ with one of positive integers $1,\ldots,p-1$ (of course, $n\ne0$). One can see that matrix~$g_3$, for a given~$n$, can be obtained by taking away from $G_3$ the rows and columns number $n$, $p$, $p+n$, $2p$, $2p+n$, $3p$, $3p+n$, $4p$, $4p+1$, and $4p+3$. Let $\widetilde S_1$ (resp.~$\widetilde S_2$) denote the $(p-2)\times (p-2)$ matrix obtained by removing the $n$th and $p$th rows and columns from the matrices $S_1$ (resp.~$S_2$). We set:
\begin{displaymath}
s_n = \det \widetilde S_1,\qquad t_n=\det \widetilde S_2.
\end{displaymath}
Also, let $\widetilde S_3$ denote the matrix obtained by removing the f\/irst row and column from the matrix~$S_3$, that is, $\widetilde S_3 = (p)$.

\begin{theorem}
\label{th_sntn} The invariant $I\bigl(L_n(p,q)\bigr)$ is explicitly given by
\begin{equation}
I\bigl(L_n(p,q)\bigr) = - \frac{1}{s_n^2 \, t_n^2 \, p^2} \label{eq_In},
\end{equation}
where
\begin{gather}
s_n = n q^*_n - p \, \nu_n, \label{eq_sn}
\\
t_n = p - s_n = p \, (\nu_n + 1) - n q^*_n, \label{eq_tn}
\\
\nu_n = \frac{1}{p} \, \sum\limits_{k = 0}^{p-1} \frac{1 - \zeta^{k(1-n)}}{1 - \zeta^k} \, \frac{1 - \zeta^{kq(q^*_n-1)}}{1 - \zeta^{-kq}}, \quad \zeta = e^{2\pi i/p}
\label{eq_nun}
\end{gather}
and $q^*_n \in \{1, \ldots, p-1\}$ is such that $q q^*_n \equiv n \bmod p$.
\end{theorem}

\begin{remark}
As we will prove, the values $\nu_n$, def\/ined in~(\ref{eq_nun}), are integers too and belong to $\{0, \ldots, n-1\}$. So, we have the congruences $s_n \equiv nq^*_n \bmod p$ and $t_n \equiv - nq^*_n \bmod p$.
\end{remark}

\begin{proof}[Proof of Theorem~\ref{th_sntn}]

Equation~(\ref{eq_In}) is directly deduced from the block structure in equation~(\ref{eq_block_structure}) of matrix $G_3$ and equation~(\ref{eq_inv0}). So, it remains only to prove the formulas \eqref{eq_sn} and~\eqref{eq_tn} for values $s_n$ and $t_n$.

We f\/irst prove~(\ref{eq_sn}). We use the factorization of matrix $S_1$ given by~(\ref{eq_S1}) in order to simplify the matrix $\widetilde S_1$ with the help of certain sequence of elementary transformations preserving the determinant.

Recall that matrix~$\widetilde S_1$ is obtained from matrix
\begin{equation}
S_1 = ({\bf 1}_p - E^q) ({\bf 1}_p - E)
\label{eq_small1}
\end{equation}
by taking away $n$th and $p$th columns and rows. This means that $\widetilde S_1$ can be obtained also as a product like~(\ref{eq_small1}),
but with the corresponding rows withdrawn from matrix~$({\bf 1}_p - E^q)$, and corresponding columns withdrawn from matrix~$({\bf 1}_p - E)$. Note that below, when we are speaking of row/column numbers in matrix~$\widetilde S_1$, we mean the numbers that these rows/columns had in~$S_1$, \emph{before} we have removed anything from it.

So, here are our elementary transformations. In matrix $({\bf 1}_p - E^q)$, for each integer $k$ from 1 to $q^*_n - 1$, we add the $(kq)$th row
to the $(kq+q)$th row (numbers modulo~$p$).
In matrix~$({\bf 1}_p - E)$, we f\/irst add the $(p-1)$th column to the $(p-2)$th one, then we add the $(p-2)$th column to the $(p-3)$th one and so forth omitting the pair of column numbers $n-1$ and $n+1$. Then, the resulting matrix has a determinant is equal to $s_n = \det \widetilde S_1$, and admits the following structure:
\begin{equation}
\left(\begin{array}{ccc} {\bf 1}_{p - 2} & {\bf c}_n & {\bf c}_p\end{array}\right) \, \left(\begin{array}{c}{\bf 1}_{p - 2} \\ {\bf r}_n \\ {\bf r}_p \end{array}\right) = {\bf 1}_{p - 2} + {\bf c}_n \otimes
{\bf r}_n + {\bf c}_p \otimes {\bf r}_p.
\label{eq_tildeS1}
\end{equation}
Here we have also moved the $n$th row in matrix $({\bf 1}_p - E^q)$ to the $(p-1)$th position and the $n$th column in matrix $({\bf 1}_p - E)$ to the $(p-1)$th position. The components of column ${\bf c}_p$ and row ${\bf r}_n$ look like
\begin{gather*}
({\bf c}_p)_i = \left\{
\begin{array}{ccc}
1, & i = kq \bmod p, & k = 1, \ldots, q^*_n - 1, \\
0, && {\rm otherwise},\hfill \\
\end{array}
\right.
\\
({\bf r}_n)_i = \left\{
\begin{array}{cc}
1, & i = 1, \ldots, n - 1, \\
0, & i = n, \ldots, p - 2. \\
\end{array}
\right.
\end{gather*}
Moreover, for all $i$ we have
\begin{equation}
({\bf c}_n)_i = 1 - ({\bf c}_p)_i, \qquad ({\bf r}_p)_i = 1 - ({\bf r}_n)_i. \label{eq_2sums}
\end{equation}

The matrix ${\bf c}_n \otimes {\bf r}_n + {\bf c}_p \otimes {\bf r}_p$ has rank 2, so its eigenvalues are
\begin{displaymath}
0, \ldots, 0, \lambda_1, \lambda_2,
\end{displaymath}
where $\lambda_1, \lambda_2$ are the eigenvalues of the following matrix of size $2\times 2$:
\begin{displaymath}
\left(
\begin{array}{cc}
{\bf r}_n {\bf c}_n & {\bf r}_n {\bf c}_p \\
{\bf r}_p {\bf c}_n & {\bf r}_p {\bf c}_p
\end{array} \right).
\end{displaymath}
Therefore, from~(\ref{eq_tildeS1}), we can deduce that the determinant of $\widetilde S_1$ is equal to
\begin{equation*}
\det \widetilde S_1 = s_n = \left|
\begin{array}{cc}
1 + {\bf r}_n {\bf c}_n & {\bf r}_n {\bf c}_p \\
{\bf r}_p {\bf c}_n & 1 + {\bf r}_p {\bf c}_p
\end{array} \right|.
\end{equation*}
Further, using~(\ref{eq_2sums}) and elementary transformations, we simplify this determinant to
\begin{displaymath}
s_n = \left|
\begin{array}{cc}
n & {\bf r}_n {\bf c}_p \\
p & q^*_n \bmod p
\end{array} \right| = n q^*_n - p \, {\bf r}_n {\bf c}_p,
\end{displaymath}
where the inner product ${\bf r}_n {\bf c}_p$ is an integer between $0$ and $n-1$.

Finally, using the discrete Fourier transform, one can prove that
\begin{displaymath}
{\bf r}_n {\bf c}_p = \nu_n = \frac{1}{p} \, \sum\limits_{k = 0}^{p-1} \frac{1 - \zeta^{k(1-n)}}{1 - \zeta^k} \, \frac{1 - \zeta^{kq(q^*_n-1)}}{1 - \zeta^{-kq}},
\end{displaymath}
where $\zeta = e^{2\pi i/p}$.

Quite similarly, we obtain formula~(\ref{eq_tn}).
\end{proof}

\subsection{Invariant for all framings}
\label{subsec_lens_framings}

According to Section~\ref{sec_framing}, we should investigate the change of matrix $G_3$ under the change of the framing. We assume that we do the
f\/irst half-revolution exactly as described in Section~\ref{sec_framing}, and the second half-revolution goes in a similar way but with the pair of
edges $AB$, $CD$ replaced by the pair~$AC$, $BD$, the third half-revolution involves again the pair~$AB$, $CD$ and so on.

Thus, we have to study how the submatrices $S_1$ and $S_2$ of $G_3$ change, because by~(\ref{eq_block_structure}) they correspond to the pairs $AB,CD$ and~$AC,BD$ respectively. We think of these matrices as made of the following blocks:
\begin{equation}
S_i = \left( \begin{array}{cc} K_i & M_i \\ L_i & \beta_i \end{array} \right), \qquad i = 1, 2, \label{eq_Si0}
\end{equation}
where $K_i$ is a $(p-1) \times (p-1)$ matrix, $L_i$ and $M_i$ are row and column of size $p-1$ respectively, and $\beta_i$ is a real number.

We keep the notations used in Section~\ref{sec_framing}. The changes made in matrices $S_1$ and $S_2$ follow from formula~(\ref{eq_***}). When we do the f\/irst half-revolution we have $\epsilon=1$ according to our agreement that the orientation of the tetrahedra in Figs.~\ref{fig1a} and~\ref{fig1b} is~$BACD$. When we do the second half-revolution we have $\epsilon=-1$, because the orientation of the ``initial'' (or better to say, the innermost in the ``sandwich'', see Section~\ref{sec_framing}) tetrahedron has changed. Then $\epsilon$ takes again the value~$-1$, and so on.

Suppose we have done this way $h$ half-revolutions. We let $S_1^{(h)}$ (resp.~$S_2^{(h)}$) denote the matrices obtained from $S_1$ (resp.~$S_2$) according to~(\ref{eq_***}). Then
\begin{equation*}
S_1^{(2m-1)} = S_1^{(2m)} = \left( \begin{array}{cccccc}
K_1 & M_1 & {} & {} & {} & {} \\ L_1 & \beta_1-1 & 1 & {} & {} & {} \\ {} & 1 & -2 & 1 & {} & {} \\ {} & {} & 1 & \ddots & \ddots & {} \\
{} & {} & {} & \ddots & -2 & 1 \\ {} & {} & {} & {} & 1 & -1 \end{array} \right), 
\end{equation*}
where the total number of $(-2)$'s is $m-1$, and
\begin{equation*}
S_2^{(2m)} = S_2^{(2m+1)} = \left( \begin{array}{cccccc}
K_2 & M_2 & {} & {} & {} & {} \\ L_2 & \beta_2+1 & -1 & {} & {} & {} \\ {} & -1 & 2 & -1 & {} & {} \\ {} & {} & -1 & \ddots & \ddots & {} \\
{} & {} & {} & \ddots & 2 & -1 \\ {} & {} & {} & {} & -1 & 1 \end{array} \right), 
\end{equation*}
where the total number of $2$'s is $m-1$. By def\/inition, $S_1^{(-1)} = S_1^{(0)} = S_1$ and $S_2^{(0)} = S_2$.

In conformity with the notations used in the previous subsection, we let $\widetilde S_1^{(2m)}$ denote the matrix obtained by taking away the
$n$th and the last columns and rows from matrix $S_1^{(2m)}$. Set $s_n^{(2m)} = \det \widetilde S_1^{(2m)}$. Quite similarly, we def\/ine $s_n^{(2m+1)} = \det \widetilde S_1^{(2m+1)}$ and $t_n^{(2m)} = \det \widetilde S_2^{(2m)}$. The following result gives the value of our invariant for all framings in~$L_n(p,q)$.
\begin{theorem}
\label{th_snmtnm}
Let $r$ be a difference between the considered and the simplest $($i.e., as in Subsection~{\rm \ref{subsec_lens_calculation})} framings of the $n$-th unknot in $L(p, q)$.

 The invariant $I^{(r)}\bigl(L_n(p,q)\bigr)$, for $r$ integer or half-integer, is given by
\begin{gather}
I^{(m)}\bigl(L_n(p,q)\bigr) = \frac{- 1}{\bigl(s_n^{(2m)}\bigr)^2 \, \bigl(t_n^{(2m)}\bigr)^2 \, p^2},
\label{eq_InG1}
\\
I^{(m + 1/2)}\bigl(L_n(p,q)\bigr) = \frac{1}{\bigl(s_n^{(2m+1)}\bigr)^2 \, \bigl(t_n^{(2m)}\bigr)^2 \, p^2},
\label{eq_InG2}
\end{gather}
where
\begin{gather}
s_n^{(2m)} = (- 1)^m (s_n - m p), \label{eq_sn2m}
\\
s_n^{(2m + 1)} = (- 1)^{m+1} (s_n - m p - p), \label{eq_sn2m1}
\\
t_n^{(2m)} = t_n + m p = p - s_n + m p. \label{eq_tn2m}
\end{gather}
\end{theorem}

\begin{proof}
By equation~(\ref{eq_inv0}), we have two formulas for the invariant:
\begin{equation*}
I^{(m)}\bigl(L_n(p,q)\bigr) = \frac{- 1}{\bigl(s_n^{(2m)}\bigr)^2 \, \bigl(t_n^{(2m)}\bigr)^2 \, p^2}
\end{equation*}
and
\begin{equation*}
I^{(m + 1/2)}\bigl(L_n(p,q)\bigr) = \frac{1}{\bigl(s_n^{(2m+1)}\bigr)^2 \, \bigl(t_n^{(2m)}\bigr)^2 \, p^2}. 
\end{equation*}
So, what remains is to specify the values of $s_n^{(2m)}$, $t_n^{(2m)}$ and~$s_n^{(2m+1)}$.

First of all, we need a lemma concerning matrices $S_1$ given by~(\ref{eq_S1}) and~$S_2$ given by~(\ref{eq_S2}). Note that they are degenerate, so they do not have inverse matrices. Instead, we can consider their \emph{adjoint} matrices, whose rank is necessarily not bigger than $1$.

\begin{lemma}
\label{lemma_small}
The adjoint matrix to both $S_1$ and $S_2$ has all its elements equal to~$p$.
\end{lemma}

\begin{proof}
One can see that the adjoint matrix to ${\bf 1}_p-E^r$, where $E$ is given by~(\ref{eq_E}) and $p$ and~$r$ are relatively prime, is a matrix whose all elements are unities. When we take a product like in~(\ref{eq_S1}) or~(\ref{eq_S2}), the corresponding adjoint matrices are also multiplied (this can be seen at once if we think of matrices $S_1$ and~$S_2$ and their factors in (\ref{eq_S1}) and~(\ref{eq_S2}) as limits of some nondegenerate matrices, keeping in mind that the adjoint to a nondegenerate matrix~$A$ is $\det A\cdot A^{-1}$). The product of two $p\times p$ matrices whose all elements are~$1$ is a matrix whose all elements are~$p$.
\end{proof}

Return to the proof of Theorem~\ref{th_snmtnm}. First, we prove~(\ref{eq_sn2m}). By def\/inition, $s_n^{(0)} = s_n$. Let us f\/ind the value
\begin{equation}
s_n^{(2)} = \det \widetilde S_1^{(2)} = \left| \begin{array}{cc} \widetilde K_1 & \widetilde M_1 \\ \widetilde L_1 & \beta_1-1 \end{array} \right|,
\label{eq_small2}
\end{equation}
where $\widetilde K_1$ means the matrix $K_1$ without its $n$th row and $n$th column; $\widetilde L_1$ and~$\widetilde M_1$ mean~$L_1$ and~$M_1$ without
their $n$th entries.

It follows from Lemma~\ref{lemma_small} that
\begin{equation}
\left| \begin{array}{cc} \widetilde K_1 & \widetilde M_1 \\ \widetilde L_1 & \beta_1 \end{array} \right| = p.
\label{eq_small3}
\end{equation}
Indeed, the matrix
\begin{displaymath}
\left( \begin{array}{cc} \widetilde K_1 & \widetilde M_1 \\ \widetilde L_1 & \beta_1 \end{array} \right)
\end{displaymath}
is just the matrix $S_1$ (compare with~(\ref{eq_Si0})) without its $n$th row and $n$th column, so the determinant~(\ref{eq_small3}) is the corresponding element of the matrix adjoint to~$S_1$. Comparing~(\ref{eq_small2}) and~(\ref{eq_small3}), we get
\begin{displaymath}
s_n^{(2)} = p - s_n.
\end{displaymath}

Let $m > 1$. Using the row (or column) expansion of the determinant of matrix $\widetilde S_1^{(2m)}$, one can see that a number sequence $s_n^{(2m)}$ is def\/ined by the recurrent condition
\begin{equation*}
s_n^{(2m)} = -2 s_n^{(2m-2)} - s_n^{(2m-4)}.
\end{equation*}
Besides, as we have just shown, there are the following initial conditions $s_n^{(0)} = s_n$, $s_n^{(2)} = p - s_n$. Hence, using induction on $m$, we get~(\ref{eq_sn2m}).

As for~(\ref{eq_sn2m1}) and~(\ref{eq_tn2m}), they are proved in much the same manner.
\end{proof}

\begin{remark}
\label{remark_L71}
If the values of our invariants~(\ref{eq_InG1}) and~(\ref{eq_InG2}) turn into inf\/inity, then this means that complex~(\ref{eq_alg_compl}) is \emph{not acyclic}. However, one can check that, e.g., for $p = 7$ this never happens. Moreover, if we consider the lens space $L(7,1)$, the invariants~(\ref{eq_InG1}) and~(\ref{eq_InG2}) are enough to distinguish \emph{all} unknots with \emph{all} framings from each other.
\end{remark}

\section{Discussion of results}
\label{sec_discussion}

Here are some remarks about the results of this paper and possible further directions of research.

{\bf Nontriviality of the invariant.} Given the results of Theorems~\ref{th_Iunknot}, \ref{th_sntn} and~\ref{th_snmtnm}, we
observe that our invariant, although being just one number, is highly nontrivial already for
the unknot in sphere~$S^3$ and ``unknots'' in lens spaces: it detects the framing of unknot
in~$S^3$ and, for example, it is powerful enough to distinguish between all ``unknots'' with
all framings in~$L(7,1)$, as stated in Remark~\ref{remark_L71}.

The work is currently underway of doing calculations for nontrivial knots in~$S^3$ and other
nontrivial situations. Regardless of the fact to what extent this invariant can detect
nontrivial knots, we have seen that it does give interesting information. Moreover, recall
that this is just the ``zeroth level invariant'' in terminology of
papers~\cite{korepanov-2009-1,korepanov-2009-2}; we leave the calculations for other levels
for further work.

Thus, although much work remains to be done, the calculations in the present paper already
demonstrate the nontriviality of our theory.

{\bf More quantum nature than in the usual Reidemeister torsion.} Our invariants in this paper, as well as in~\cite{korepanov-2009-1,korepanov-2009-2},
are calculated using the torsion of acyclic complexes of geometric origin. Comparing this
with the ``usual'' Reidemeister torsion, we see a striking dif\/ference between them. Our
invariants work like quantum invariants in the sense that they do not require any nontrivial
representation of the fundamental group (be it the fundamental group of the manifold or of
the knot complement). We expect that a ``non-commutative'' version of our invariants will be
developed, where real or complex numbers attached to elements of a~triangulation are
replaced with some associative algebra, which will show even more clearly the quantum nature
of our theory.

{\bf Using representations of the fundamental group.} There also exists a version of our theory using the universal cover of a manifold~$M$ and nontrivial representations of the fundamental group $\pi_1(M)$ into the group of motions
of three-dimensional Euclidean space~\cite{KM,M,M-thesis}. This theory has been developed
yet only for closed manifolds. In contrast with the previous remark, the invariants in such
theory appear to be related to the usual Reidemeister torsion, including even the
non-abelian case, see~\cite{M}; it is also interesting to compare this with
paper~\cite{JDFourier}. As stated in~\cite{JDFourier}, the Reidemeister torsion is related
to the \emph{volume form} on the character variety of the fundamental group. Note also that
(see for example~\cite{K}) the structure of the character variety of the fundamental group
of the manifold obtained by a surgery on a knot is well-known: it can be deduced using the
character variety of the knot complement and the slope of the surgery. So, we can try to
search for possible links between that and the set of real numbers we obtain using our
invariant, guided by a general idea that, in a sense, the representations could be hidden in
the surgery.

{\bf Framings, surgeries and shellings.} One more direction of research is suggested by the presence of a framed knot in our
constructions. As we know, this is usually used for obtaining new closed manifolds by means
of a surgery. For example, lens spaces are obtained by surgery on the unknot in~$S^3$ with a
certain framing, which is exactly the slope of the surgery. So, the idea is, in a general
formulation, to explore more in-depth the behavior of our invariant under surgeries. This
can require more research on what happens with our invariants under shellings of a manifold
boundary.

Shellings and inverse shellings of a manifold boundary should be studied also because they
provide a way of changing the boundary triangulation,
compare~\cite[Section~4]{korepanov-2009-1}.

{\bf Other geometries, TQFT's and higher-dimensional manifolds.} Finally, we would like to emphasize once more that, from a purely mathematical point of view, three-dimensional Euclidean geometry plays in this paper just an auxiliary role as a provider of algebraic structure~-- a chain complex made of dif\/ferentials of Euclidean quantities. In general, this auxiliary geometry need not be Euclidean, nor three-dimensional, as is shown by other theories developed for instance in~\cite{SL2,3dim4geo,KKM,BKM}. Note that the theory in paper~\cite{BKM} has already been developed to the stage of a working TQFT.

Especially interesting must be generalization to \emph{four}-dimensional manifolds. The work in this direction is currently underway~\cite{K3324}.

\subsection*{Acknowledgements}

\looseness=1
The authors J.D.\ and I.G.K.\ acknowledge support from the \emph{Swiss National Science
Foundation}. In particular, it made possible the visit of I.G.K.\ to Geneva, where the work on
this paper began.
The work of I.G.K.\ and E.V.M.\ was supported by \emph{Russian Foundation
for Basic Research}, Grants no.~07-01-96005 and 10-01-96010.
The work of J.D. was partially supported by the {European Community} with Marie Curie
Intra--European Fellowship (MEIF--CT--2006--025316). While writing the paper, J.D.\ visited
the CRM. He thanks the CRM for hospi\-tality.

Finally, we would like to thank the referee for the very attentive reading of the manuscript and constructive criticism. We hope that, following the referee's comments, we have substantially improved our paper.

\pdfbookmark[1]{References}{ref}
\LastPageEnding

\end{document}